\documentclass[preprint]{elsarticle}

\usepackage{lipsum}
\makeatletter
\def\ps@pprintTitle{%
 \let\@oddhead\@empty
 \let\@evenhead\@empty
 \def\@oddfoot{}%
 \let\@evenfoot\@oddfoot}
\makeatother

\usepackage{graphicx}
\usepackage{amsmath}
\usepackage{amssymb}
\usepackage{enumerate}
\usepackage{float}
\usepackage{color}

\usepackage{amsthm}
\usepackage[T1]{fontenc}
\usepackage[para]{threeparttable}
\usepackage[margin=1in]{geometry}
\usepackage{booktabs}
\usepackage{multirow}
\usepackage{longtable}
\usepackage{booktabs}
\usepackage{tabularx}
\usepackage{lscape}
\usepackage{threeparttable}
\usepackage{tablefootnote}
\usepackage{tikz}
\usepackage{amsmath}
\usepackage{algorithm}
\usepackage[noend]{algpseudocode}
\usetikzlibrary{positioning}
\usepackage{hyperref}
\makeatletter
\usepackage[labelformat=simple]{subcaption}

\usepackage{multirow}
\usepackage{lscape}
\usepackage{longtable}
\usepackage{lineno,hyperref}
\usepackage{setspace}
\modulolinenumbers[5]
\usepackage{endnotes}
\let\footnote=\endnote

\usepackage{natbib}
 \bibpunct[, ]{(}{)}{,}{a}{}{,}%

\usepackage{array}

\usepackage{pdflscape}

\def\comp{\raise 1pt \hbox{$\scriptstyle\circ$}}

\def\argmax{\mathop{\rm argmax}\limits}

\def\min{\mathop{\rm min\ }\limits}
\def\max{\mathop{\rm max\ }\limits}
\def\st{\mathop{\rm subject\ to\ }}

\def\upto{{\raise 1pt \hbox{$\scriptstyle \,\nearrow\,$}}}
\def\downto{{\raise 1pt \hbox{$\scriptstyle \,\searrow\,$}}}

\newtheorem{theorem}{Theorem}

\biboptions{authoryear}	

\begin{document}
\begin{frontmatter}
	\title{An Exact Algorithm for Optimization Problems with Inverse S-shaped Function}
 \author{Arka Das}
 \ead{phd17arkadas@iima.ac.in}
\author{Ankur Sinha\corref{cor1}}
\ead{asinha@iima.ac.in}
\author{Sachin Jayaswal}
\ead{sachin@iima.ac.in}
\cortext[cor1]{Corresponding author}
\address{Operations and Decision Sciences, Indian Institute of Management Ahmedabad, Gujarat, India, 380015}
\begin{abstract}

In this paper, we propose an exact general algorithm for solving non-convex optimization problems, where the non-convexity arises due to the presence of an inverse S-shaped function. The proposed method involves iteratively approximating the inverse S-shaped function through piece-wise linear inner and outer approximations. In particular, the concave part of the inverse S-shaped function is inner-approximated through an auxiliary linear program, resulting in a bilevel program, which is reduced to a single level using KKT conditions before solving it using the cutting plane technique.
To test the computational efficiency of the algorithm, we solve a facility location problem involving economies and dis-economies of scale for each of the facilities. The computational experiments indicate that our proposed algorithm significantly outperforms the previously reported methods. We solve non-convex facility location problems with sizes up to 30 potential facilities and 150 customers. The computational experiments indicate that our proposed algorithm converges to the global optimum within a maximum computational time of 3 hours for 95$\%$ of the datasets. For almost 60 percent of the test cases, the proposed algorithm outperforms the benchmark methods by an order of magnitude. The paper ends with managerial insights on facility network design involving economies and dis-economies of scale. One of the important insights points out that it may be optimal to increase the number of production facilities operating under dis-economies of scale with an overall decrease in transportation costs.
\end{abstract}
\begin{keyword}
Facility Location, Inverse S-shaped Function, Economies and Dis-economies of Scale
\end{keyword}
\end{frontmatter}
\section{Introduction} \label{Intro}
In this paper, we study the general non-convex problems, where non-convexity arises due to the presence of inverse S-shaped function(s). Here, we consider an optimization model of the following kind:
	\begin{align}
	\min_{x,y,z} & \; f(x,y) + \sum_{j=1}^{J}\phi_j(z_j) \label{eq:origObj}\\
	\st & g_{m}(x,y,z)\leq 0, \quad  m=1,\dots,M; \label{eq:origConst}\\
    & e_j \leq z_j \leq K_j, \qquad j=1,\dots,J \label{eq:origConst2}\\
    & x \in  \mathbb{R}^{I \times J} \label{eq:origConst3}\\
    & z \in  \mathbb{R}^{ J} \label{eq:origConst4}\\
	& y \in \mathbb{Z}^n \label{eq:origConst5}
	\end{align}
where $f(x,y)$ and $g_{m}(x,y,z)$ are convex functions, $\phi_j(z_j)$ is an inverse S-shaped function, and constraint \ref{eq:origConst2} sets the bounds for the variable $z_j$. In \eqref{eq:origConst3} - \eqref{eq:origConst5} we assume variables $x$ and $z$ to be continuous, and $y$ to be integers. In this problem, the source of non-convexity arises from the inverse S-shaped function in the objective. The inverse S-shaped function is neither fully concave nor fully convex; rather, it is concave in nature up to a certain point and then transitions into convex at a point that is often referred to as the deflection point or economic point. \par

Non-convex problems have been of interest for the last sixty years due to the presence of several application problems that contain various kinds of non-convexities. These problems may contain multiple local optimums because of which they are often difficult to solve to global optimality and therefore have kept both researchers and practitioners busy. Some of the practical problems that are non-convex include, facility location problems with economies and dis-economies of scale \citep{lu2014facility}, location-inventory network design problems \citep{li2021location,farahani2015location,shen2007incorporating,jeet2009logistics}, facility location problems with concave cost \citep{soland1974optimal,ni2021branch}, production transportation problems \citep{holmberg1999production,kuno2000lagrangian,saif2016supply}, Stochastic Transportation-Inventory Network Design Problem \citep{shu2005stochastic}, concave minimum cost network flow problems \citep{guisewite1990minimum,fontes2007heuristic}, etc. \par
Researchers have reported several algorithms and heuristics for solving non-convex problems, but most of them do not scale well while solving high-dimensional instances. The following methods are some of the exact algorithms that have been reported to solve non-convex problems: piece-wise linear underestimation \citep{dasci2001plant}, branch and bound \citep{falk1969algorithm,soland1971algorithm, shectman1998finite}, spatial branch and bound \citep{mccormick1976computability,lee2001global,belotti2009branching}, branch and reduce \citep{ryoo1995global,ryoo1996branch}, cutting plane algorithm \citep{tuy1964concave,taha1973concave,tawarmalani2005polyhedral}, hybrid method \citep{marsten1978hybrid}, reformulation-linearization based technique \citep{zetina2021exact}, etc. Most of these algorithms are suited to solve either concave minimization problems or very specific kinds of non-convex problems. To the best of our knowledge, these algorithms are not suitable to solve the general class of non-convex problems where non-convexity arises due to inverse S-shaped function(s). Problems with S-shaped functions appear in various kinds of application problems, like facility location problems \citep{lu2014facility}, location inventory problems \citep{li2021location}, knapsack problems \citep{augrali2009solving}, etc. \cite{lu2014facility} was the first to report a Lagrangian relaxation-based column generation heuristic for solving facility location problems with inverse S-shaped functions. Though the method was able to find good solutions, the method is neither exact nor general in nature. Recently, \cite{li2021location} assumed warehouse operating cost as an inverse S-shaped function of the volume of demand assigned to a warehouse and modeled it as a location and two-echelon inventory network design problem. They formulated the problem as a set covering problem and solved its linear relaxation using a column generation algorithm. Though they claim to solve the pricing problem of the column generation using a branch-and-bound method by exploiting its structural properties, they did not report the results after solving the pricing problem. The method is not suited to handle problems with general inverse S-shaped functions; for instance, it cannot handle non-differentiable inverse S-shaped functions. 

%
In this paper, we propose an exact algorithm to solve large-scale instances of problems with general forms of inverse S-shaped cost function(s). The method can be used to solve all classes of non-convex problems arising in the area of operations management, where non-convexity often arises due to the presence of economies and dis-economies of scale. The method converts the original non-convex problem to a relaxed mixed integer convex formulation by replacing the concave part of the inverse S-shaped function with its inner approximation using certain principles of bilevel programming. The relaxed convex formulation provides a lower bound to the original problem when solved exactly using the cutting plane algorithm (discussed in Appendix C).
The above steps are repeated iteratively, guaranteeing improvement of the lower bound at each iteration and convergence to the global optimum of the original problem upon termination. The methodology is tested on facility location problems of dimensions up to 30 facilities and 150 customers. We compare the computational performance of the proposed method against \cite{lu2014facility}. The proposed general method outperforms the benchmark that exploits problem specific information by an order of magnitude on 60$\%$ of the test instances. The paper is structured as follows.
In Section~\ref{sec:SolutionMethod} we discuss the algorithm followed by the application problem in Section~\ref{sec:Appprob}. Section~\ref{sec:CompExper} contains the data synthesis scheme and computational results. The paper also provides useful managerial insights on the design of facility networks in the presence of economies and dis-economies of scale. Finally, we conclude our findings in Section~\ref{sec:conclusions}. We also illustrate our proposed method through an example problem in Appendix B.

\section{Solution Method}\label{sec:SolutionMethod}
In this section, we explain the working of our proposed algorithm that solves \eqref{eq:origObj}-
\eqref{eq:origConst5}. As explained earlier, this problem is non-convex in nature due to the presence of $\phi_j(z_j)$. The initial shape of $\phi_j(z_j)$ is strictly concave up to a point ($z_j^0$), after which the shape turns strictly convex. The point $z_j^0$ is referred to as the deflection point/economic point. The function $\phi_j(z_j)$ is continuous at $z_j^0$ but not necessarily differentiable. $K_j$ is the upper bound for $z_j$. 
	\begin{figure}
		   \centering
		\begin{subfigure}{.45\textwidth}
			\centering
			\includegraphics[width=0.9\linewidth]{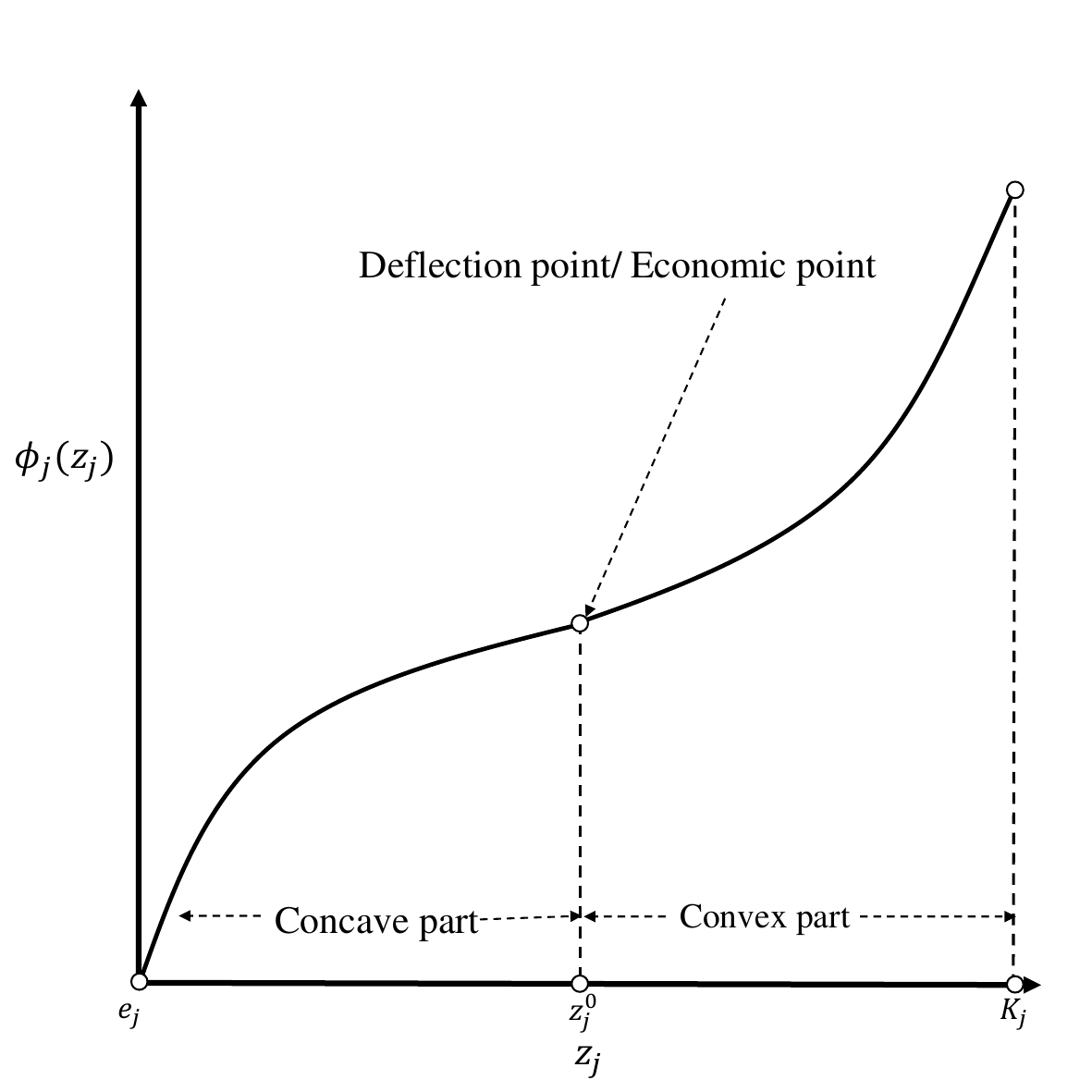}
			\caption{Inverse S-shaped function, $\phi_j(z_j)$}
			\label{fig:inverse_S_fun}
		\end{subfigure}%
		\begin{subfigure}{.5\textwidth}
			\centering
			\includegraphics[width=1\linewidth]{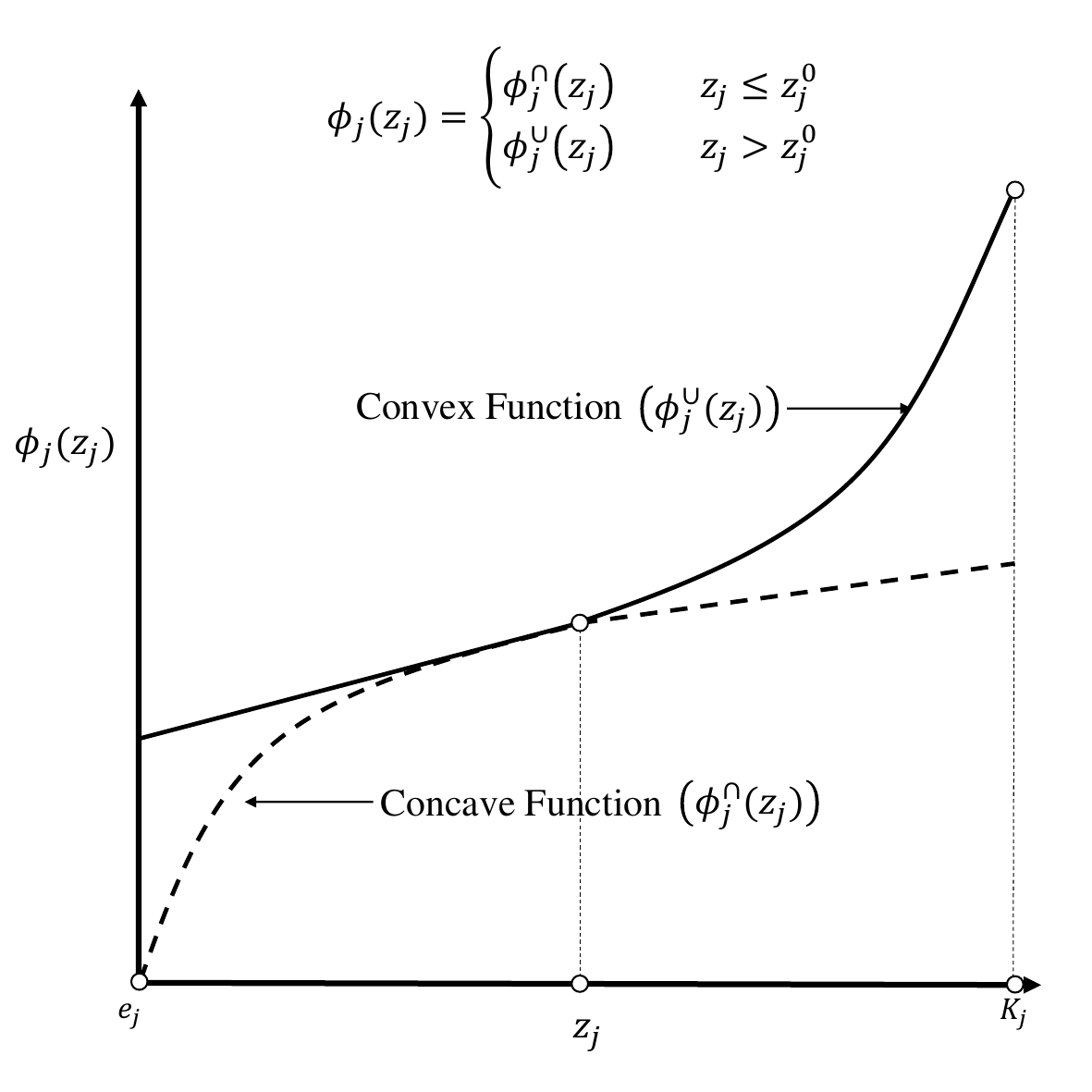}
			\caption{$\phi_{j}(z_j)$ separated into two functions $\phi_{j}^{\cap}(z_j)$ and $\phi_{j}^{\cup}(z_j)$}
			\label{fig:sep_s_fun}
		\end{subfigure}
		\caption{ Process of converting the inverse S-shape function into separate concave and convex function}
		\label{fig:convert_s_func}
	\end{figure}
		Figure~\ref{fig:inverse_S_fun} shows a function $\phi_j(z_j)$ that is continuous but non-differentiable at $z_j^0$. Note that the function ($\phi_j(z_j)$) may be non-differentiable at other points as well but for simplicity we have shown an inverse S-shaped function that is non-differentiable only at $z_j^0$. Such a function can be written in terms of a concave function and a convex function, $\phi_j^{\cap}(z_j)$ and $\phi_j^{\cup}(z_j)$, respectively, as shown in Figure~\ref{fig:sep_s_fun}. 
		If the function $\phi_j(z_j)$ is not differentiable at $z_j^0$, we need to add different tangents on the left hand side and right hand side of $z_j^0$. Let $\phi_{j_{-}}^{'}(z_j^0)$ and $\phi_{j_{+}}^{'}(z_j^0)$ be the left hand derivative and the right hand derivative of the function $\phi_j(z_j)$ at $z_j^0$. Then $\phi_{j}^{\cap}(z_j)$ and $\phi_{j}^{\cup}(z_j)$ can be defined as follows.\\
	\begin{equation}
	\phi_{j}^{\cap}(z_j) = 
	\begin{cases}
	\phi_j(z_j) & \forall \; z_j \in [e_j,z_j^0]\\
	\phi_j(z_j^0) + (z_j-z_j^0) \phi_{j_{-}}^{'}(z_j^0) \ & \forall \; z_j \in [z_j^0,K_j] \label{eq:fundef1}
	\end{cases}       
	\end{equation}
	\begin{equation}
\phi_{j}^{\cup}(z_j,z_j) =
	\begin{cases}
	\phi_j(z_j^0) + (z_j-z_j^0) \phi_{j_{+}}^{'}(z_j^0) & \forall \; z_j \in [e_j,z_j^0]\\
	\phi_j(z_j) & \forall \; z_j \in [z_j^0,K_j]\label{eq:fundef2}
	\end{cases}       
	\end{equation}
Then formulation \eqref{eq:origObj}-\eqref{eq:origConst5} can be rewritten with the help of updated $\phi_j^{\cap}(z_j)$ and $\phi_j^{\cup}(z_j)$ as:
	\begin{align}
	\min_{x,y,z,l} & \; f(x,y) + \sum_{j=1}^{J} l_{j0} \phi_{j}^{\cap}(z_j) + \sum_{j=1}^{J}l_{j1} \phi_{j}^{\cup}(z_j) \label{eq:origObj2}\\
	\st & \eqref{eq:origConst} - \;\eqref{eq:origConst5} \notag\\
	&l_{j0}+l_{j1} \leq 1, \; j=1,\dots,J \label{eq:origConst6}  \\
	& z_j^0l_{j1}\leq z_j \leq K_jl_{j1}+z_j^0 l_{j0},\; j=1,\dots,J \label{eq:origConst7}  \\
	& l_{j0},l_{j1} \in \{0,1\} ,\; j=1,\dots,J \label{eq:origConst8}
	\end{align}
	$l_{j0} \phi_{j}^{\cap}(z_j)$ and $ l_{j1} \phi_{j}^{\cup}(z_j)$ are non-convex terms, each consisting of a product of a binary variable with a continuous term. The above model can be reformulated after linearization of these terms as follows.
	\begin{align}
	\min_{x,y,z,l,w,p} & \; f(x,y) + \sum_{j=1}^{J}w_j + \sum_{j=1}^{J}p_j   \label{eq:origObj3}\\
	\st & \eqref{eq:origConst} - \;\eqref{eq:origConst5},\eqref{eq:origConst6} - \;\eqref{eq:origConst8} \label{eq:origConst9.0}\\
	& M_j^0l_{j0} \geq w_j,\; j=1,\dots,J \label{eq:origConst9}\\
	& M_j^0l_{j0}+ \phi_{j}^{\cap}(z_j) - M^0_j \leq w_j, \; j=1,\dots,J \label{eq:origConst10}\\
	& M_j^1l_{j1} \geq p_j, \; j=1,\dots,J \label{eq:origConst11}\\
	& M_j^1 l_{j1} + \phi_{j}^{\cup}(z_j) - M^1_j \leq p_j, \; j=1,\dots,J \label{eq:origConst12}
	\end{align}
The above formulation is a mixed integer program that contains separate concave ($\phi_{j}^{\cap}(z_j)$) and convex terms ($\phi_{j}^{\cup}(z_j)$). The non-convexity in the formulation (\eqref{eq:origObj3}-
\eqref{eq:origConst12}) arises because of the presence of concave terms.
We intend to solve this formulation, for which we first provide tight values for $M_j^0$ and $M_j^1$ in Appendix A, and then we discuss the handling of the concave terms in Section~\ref{sec:algorithmDescription}.
\subsection{Algorithm Description}\label{sec:algorithmDescription}
\begin{figure}
		\centering
			\includegraphics[width=0.6\linewidth]{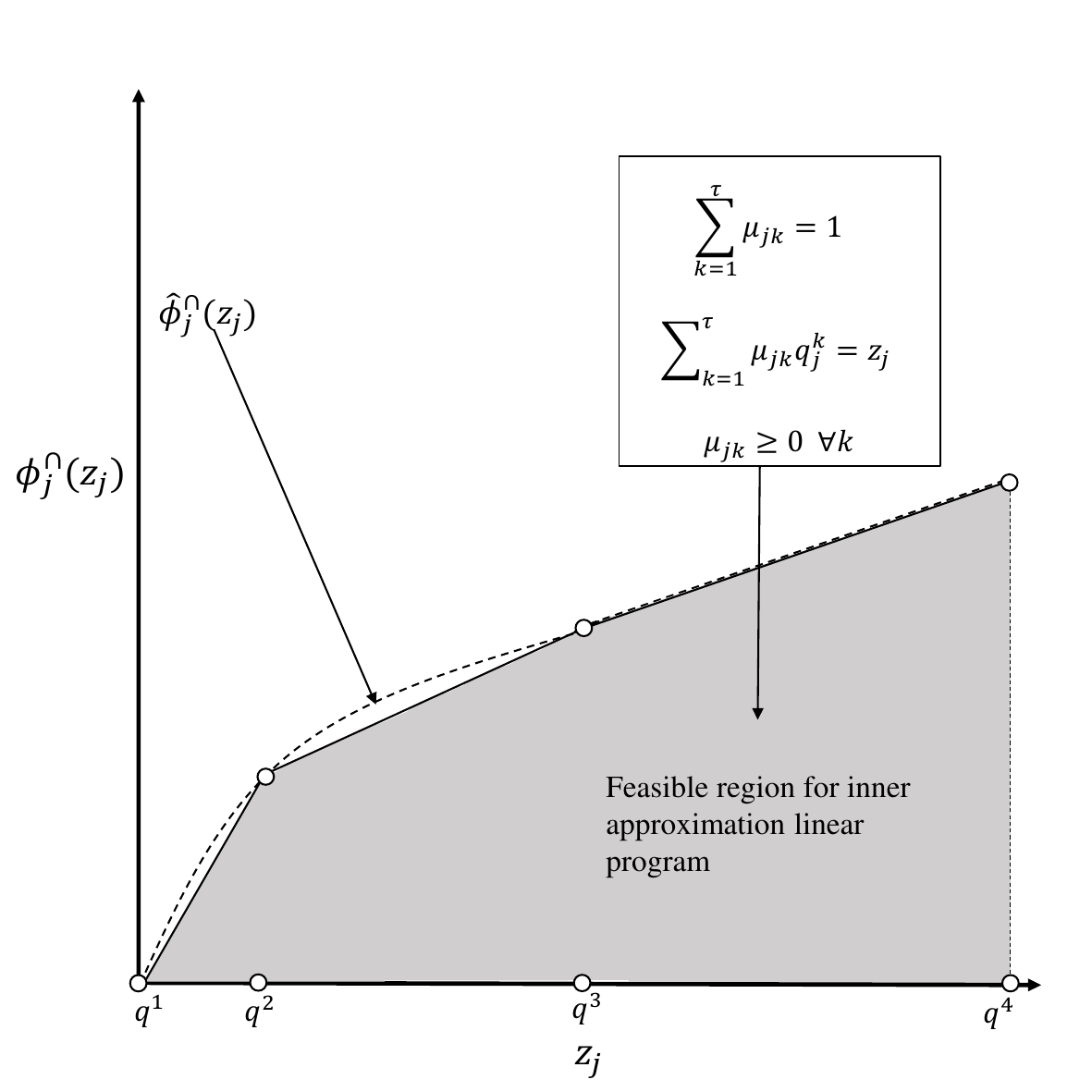}
			\caption{Inner approximation of $ \phi_j^{\cap}$ with a set of points $(\tau =4 )$}
			\label{fig:innerapproximation4points}
	\end{figure}

In this section, we discuss a solution approach for the formulation \eqref{eq:origObj3}-
\eqref{eq:origConst12} that is equivalent to \eqref{eq:origObj}-\eqref{eq:origConst5}.
Each function $\phi_j^{\cap}(z_j)\;(\forall j \in \{1,\dots, J\})$ can be inner-approximated \citep{rockafellar1970convex} using a given set of points $S_j^c=\{ q_j^{1}, q_j^{2}, \ldots, q_j^{\tau} \}$ (let $c=1$):
\begin{align}
\hat{\phi}_j^{\cap}(z_j|S_j^c) = \max_{\mu_j} \left\{ \sum_{k=1}^{\tau} \mu_{jk} \phi_{j}^{\cap}(q_{j}^{k}): \sum_{k=1}^{\tau} \mu_{jk}=1,  \sum_{k=1}^{\tau} \mu_{jk} q_j^{k}=z_{j}, \mu_{jk}\ge0, k=1,\ldots,\tau \right\}, \;\;\; \forall j \in \{1,\dots, J\} \label{eq:innerapproximation}
\end{align}
which is an inner approximation linear program with $z_j$ as a parameter and $\mu_j$ as a decision vector. For brevity, we represent the approximation $\hat{\phi}_j^{\cap}(z_j|S_j^c)$ as $\hat{\phi}_j^{\cap}(z_j)$. Figure~\ref{fig:innerapproximation4points} shows graphically the approximation of a concave function with 4 points in the approximation set using the above approach. 
The approximation \eqref{eq:innerapproximation} converts the non-convex minimization problem (\eqref{eq:origObj3}-\eqref{eq:origConst12}) into the following lower bound convex program, since $\hat{\phi}_j^{\cap}(z_j)$ is the piece-wise linear inner approximation of $\phi_j^{\cap}(z_j)$.
\begin{align}
	\min_{x,y,z,l,w,p,\zeta} & \; f(x,y) + \sum_{j=1}^{J}w_j + \sum_{j=1}^{J}p_j  \label{eq:origObj4} \\
	\st & \eqref{eq:origConst} - \;\eqref{eq:origConst5},\eqref{eq:origConst6} - \;\eqref{eq:origConst8},\eqref{eq:origConst9}, \eqref{eq:origConst11} - \;\eqref{eq:origConst12} \notag\\
	& M_j^0l_{j0}+ \zeta_j - M^0_j \leq w_j, \; j=1,\dots,J \label{eq:origConst13}\\
	& \zeta_j = \max_{\mu_j} \left\{ \sum_{k=1}^{\tau} \mu_{jk} \phi_{j}^{\cap}(q_{j}^{k}): \sum_{k=1}^{\tau} \mu_{jk}=1,  \sum_{k=1}^{\tau} \mu_{jk} q_j^{k}=z_{j}, \mu_{jk}\ge0, k=1,\ldots,\tau \right\}, \; j=1,\dots,J \label{eq:origConst14}
\end{align}
	The above lower bound formulation for \eqref{eq:origObj3}-\eqref{eq:origConst12} is a bilevel program \citep{sinha2017review} with \eqref{eq:origConst14} representing the lower level program.
	By substituting \eqref{eq:origConst14} with KKT conditions, the bi-level program is converted into single level program, where $\alpha_j$ is the Lagrange multiplier for $\sum_{k=1}^{\tau} \mu_{jk}=1$, $\beta_j$ is the Lagrange multiplier for $\sum_{k=1}^{\tau} \mu_{jk} q_j^{k}=z_{j}$, and $\gamma_{jk}$ is the Lagrange multiplier for  $\mu_{jk} \geq 0$. Thereafter, the above program reduces to the following:
	\begin{align}
	\textbf{Mod-$S^c$} \quad \min_{x,y,z,l,w,p,\zeta,\alpha,\beta,\gamma,\mu} & \;  f(x,y) + \sum_{j=1}^{J}w_j + \sum_{j=1}^{J}p_j   \label{eq:origObj5}\\
	\st &  \eqref{eq:origConst} - \;\eqref{eq:origConst5},\eqref{eq:origConst6} - \;\eqref{eq:origConst8},\eqref{eq:origConst9}, \eqref{eq:origConst11} - \;\eqref{eq:origConst12},\eqref{eq:origConst13} \label{eq:origConst15.0} \\
	& \zeta_j \geq  \sum_{k=1}^{\tau} \mu_{jk} \phi_{j}^{\cap}(q_{j}^{k}) && \hspace{-10mm} j=1,\dots,J \label{eq:origConst15} \\
	& \sum_{k=1}^{\tau} \mu_{jk}=1 && \hspace{-10mm} j=1,\dots,J  \label{eq:origConst16}\\
	& \sum_{k=1}^{\tau} \mu_{jk} q_j^{k}=z_{j} && \hspace{-10mm} j=1,\dots,J  \label{eq:origConst17}\\
	& \phi_{j}^{\cap}(q_j^{k}) - \alpha_j -  \beta_{j} q_j^k + \gamma_{jk} = 0  && \hspace{-10mm}  k=1,\ldots,\tau;  \;j=1,\ldots,J \label{eq:origConst18} \\
	& \mu_{jk}\gamma_{jk}=0 && \hspace{-10mm} k=1,\ldots,\tau;  \;j=1,\ldots,J\label{eq:origConst19}\\
	& l_{j0}, l_{j1} \in \{0,1\} && \hspace{-10mm} j=1,\ldots,J \label{eq:origConst20}\\
	& \gamma_{jk} \geq 0 && \hspace{-10mm}  k=1,\ldots ,\tau ;\;j=1,\ldots,J  \label{eq:origConst21}\\
	& \mu_{jk} \ge 0 && \hspace{-10mm} k=1,\ldots,\tau; \;j=1,\ldots,J \label{eq:origConst22}
\end{align}
	The above model consists of product terms in the complementary slackness conditions (\eqref{eq:origConst19}), which we linearize using binary variables ($u_{jk}$) and bigM ($M_{j}^2$ and $M_{j}^3$) as shown in \eqref{eq:lin1}-\eqref{eq:lin3}.
	\begin{align}
    & \quad\quad\quad \quad\quad\gamma_{jk} \leq  M_j^2 u_{jk} &&k=1,\ldots,\tau; \;j=1,\ldots,J \label{eq:lin1}\\
	& \quad\quad\quad\quad\quad \mu_{jk} \leq M_j^3 (1-u_{jk}) &&k=1,\ldots,\tau; \;j=1,\ldots,J \label{eq:lin2}\\
	& \quad\quad\quad\quad\quad u_{jk} \in \{0,1\}&& k=1,\ldots,\tau; \;j=1,\ldots,J \label{eq:lin3}
\end{align}
From constraints \eqref{eq:origConst16} and \eqref{eq:origConst22}, it is observed that the maximum value taken by $\mu_{jk}$ is $1$. Hence, $M_j^3=1$ is acceptable. $M_j^2$ can be any big number specific to the problem.

After the linearization of complementary slackness conditions, the problem \eqref{eq:origObj5}-\eqref{eq:origConst18}, \eqref{eq:origConst20} -\eqref{eq:lin3} represents a convex mixed-integer program (CMIP). 
 We solve the CMIP using cutting plane technique \citep{kelley1960cutting,jayaswal2017cutting} thereby generating an optimal solution, which is also a feasible solution to \eqref{eq:origObj}-\eqref{eq:origConst5}. 
 While the optimal value of the formulation \eqref{eq:origObj5}-\eqref{eq:origConst18}, \eqref{eq:origConst20} -\eqref{eq:lin3} provides a lower bound, substitution of the optimal solution into the objective function of \eqref{eq:origObj}-\eqref{eq:origConst5} generates an upper bound to the original problem. The optimal solution of \eqref{eq:origObj5}-\eqref{eq:origConst18}, \eqref{eq:origConst20} -\eqref{eq:lin3} generates a new approximation point $ q_j^{\tau +1}$ that is incorporated into the approximation set $S_{c}$, thereby creating a new set, $S_j^{c+1}=S_j^{c} \cup q_j^{\tau +1}$. Thereafter, $S_{c+1}$ is utilized to formulate a new CMIP, which when solved is expected to provide an improved lower bound. We keep adding new solutions to the approximation set iteratively for generating improved bounds until the optimality gap is closed. This algorithm is referred to as the Inner-Outer Approximation (IOA) algorithm in the rest of the paper. A psuedo-code of IOA algorithm is provided in Algorithm~\ref{IOAalgo}.\par
\begin{algorithm}
	\caption{IOA Algorithm for solving inverse S-shape non-convex problem}\label{IOAalgo}
	\begin{algorithmic}[1]
		\State \textit{Begin}
		\State $UB_\mathcal{A} \gets +\infty, LB_\mathcal{A} \gets -\infty, c \gets \textit{1}$
		\State $ \text{Choose an initial set of $\tau$ points } S^c = \{ q^{1}, q^{2}, \ldots, q^{\tau} \} $
		\While {($UB_\mathcal{A}-LB_\mathcal{A})/LB_\mathcal{A} > \epsilon \;  $ begin} 
		\State Solve Mod-$S^c$ (\eqref{eq:origObj5}-\eqref{eq:origConst22}) using the cutting plane technique after linearizing \eqref{eq:origConst19}
		\State Let the optimal solution for Mod-$S^c$ be $x^{\tau +c},y^{\tau +c},q^{\tau+c},w^{\tau +c},p^{\tau +c}\newline \text{\qquad for the variables $x,y,z,w,p$, respectively}$
		\State $LB_A \gets f(x^{\tau +c},y^{\tau +c}) + \sum_{j=1}^{J}w^{\tau +c}_j + \sum_{j=1}^{J}p^{\tau +c}_j $
		\State $ UB_A \gets f(x^{\tau +c}, y^{\tau +t}) + \sum_{j=1}^{J}\phi_j(q_j^{\tau +c}) $
		\State $ \mathcal{S}^{c+1} \gets  \mathcal{S}^{c} \; \cup \;q^{\tau +c}  $
		\State $c \gets c +1$
		\EndWhile
		\State \textit{End}
	\end{algorithmic}
\end{algorithm}
The algorithm starts with an initial set of points $S_j^1 =\{ q_j^{1}, q_j^{2}, \ldots, q_j^{\tau} \}$, such that $q_j^1=e_j$ and $q_j^{\tau}=K_j$, while other points may be intermediate values between $e_j$ and $K_j$. This essentially indicates that we need to pick corner points of the constraint $e_j \leq z_j \leq K_j $ to initialize the approximation set, and additional points may include any randomly drawn points between the upper bound and the lower bound. 

Note that the method guarantees global convergence, as a new point at every iteration leads to an improvement in the inner-approximation of $\phi_{j}^{\cap}(.)$ thereby leading to a tighter lower bound. In case a new point is not generated in the succeeding iteration, or two similar points are generated in two consecutive iterations, it will lead to the closure of the lower and upper bound resulting in the termination of the algorithm. We provide a theorem and proof for this claim in Appendix A.


\section{Facility location problem with economies and dis-economies of scale} \label{sec:Appprob}
The facility location problem (FLP) is a classical problem, which entails determining the locations of facilities and assigning the customers to them so as to meet the demands of the customers at a minimum total cost that constitutes of transportation costs between the facilities and their assigned customers, fixed facility opening costs, and facility operating costs. 
 To the best of our knowledge, \cite{feldman1966warehouse} was the first to consider facility fixed cost and facility operating cost to be a concave function of throughput and solved this problem using heuristics. They asserted that the concavity arises since the marginal cost of building and operating the facilities tends to decrease with an increase in throughput volume due to economies of scale. Several more recent papers in the literature have captured this behavior using a concave function to model the variable cost \citep{soland1974optimal, hajiaghayi2003facility}. \cite{hajiaghayi2003facility} introduced a general version of facility location problems by assuming facility cost (cost to connect at least one client to the facility) to be a concave function of the number of clients assigned to a facility. Several different cost functions become concave due to economies of scale, including, concave production cost \citep{cohen1991integrated,holmberg1999production,kuno2000lagrangian,dasci2001plant,romeijn2010integrating,sharkey2011exact,saif2016lagrangian}, concave inventory cost \citep{daskin2002inventory}, concave transportation cost \citep{lin2006strategic}, concave operating cost \citep{hajiaghayi2003facility}, concave site dependant costs \citep{dupont2008branch}, and concave technology acquisition cost \citep{dasci2001plant}. However, many of the production technologies also exhibit dis-economies of scale due to the overuse of resources \citep{lu2014facility}, warehouse congestions \citep{desrochers1995congested}, increasing operational complexity \citep{mcafee1995organizational}, and fixed managerial ability \citep{alvarez2003diseconomies}. Dis-economies of scale are usually modeled in the literature using a convex function. For example, congestion at facilities are often modeled using convex functions in location models \citep{desrochers1995congested,amiri1997solution,elhedhli2006service,vidyarthi2014efficient,bhatt2021alternate}. Additionally, the convex location model also appears in stochastic transportation problems \citep{holmberg1995efficient}, strategic inventory-location problems\citep{benjaafar2004demand} etc. A stream of literature also uses functions that are neither concave nor convex to model the non-linear behavior of the variable costs. For example, \cite{holmberg1994solving} modeled this behavior using stair-wise function facility location problems, while \cite{elhedhli2018service} modeled a service system design problem with stochastic demand where such non-convex functions typically arise due to the stochastic nature of the demand of the customers and variability in service rate. \par
The inverse S-shaped function is another kind of function, which is neither concave nor convex. The concept of the S-shaped cost function was first asserted by \cite{mirchandani1989discrete}. \cite{mirchandani1989discrete} modeled discrete facility location problems and briefly discussed the relationship between the facility setup cost and the number of facilities. Subsequently, \cite{lu2014facility} asserted that in a large production house, the marginal cost often starts increasing after a certain production volume. Hence initially, the production cost is concave in nature due to economies of scale, but after a specific production volume, the system starts over-utilizing the resources resulting in a convex behavior due to dis-economies of scale. As a result, the production cost function behaves like an inverse S-shaped function that is neither fully concave nor fully convex. \cite{lu2014facility} formulated the capacitated facility location problem as a mixed integer non-convex formulation, which was further decomposed into separate concave and convex subproblems using Lagrangian relaxation. The subproblems were first solved using heuristics, and then different column generation techniques were used to find the best feasible solution. Though the method performed satisfactorily in finding quality solutions, the method is neither exact nor general in nature. Hence, the method reported large gaps in termination for most of the data instances. In this paper, we are interested in solving exactly the facility location problem with economies and dis-economies of scale similar to \cite{lu2014facility}. The cost function includes production costs of the facilities, fixed setup costs of the facilities, and transportation costs between the facilities and customers. Here, the facility's production cost function behaves like an inverse S-shaped function due to economies and dis-economies of scale. \par 
To formally define the problem, let M denote a set of $m$ locations where the facility can be opened, and $N$ denote a set of $n$ customer locations. Let $c_{ij}$ be the transportation cost of satisfying demand of customer $i\in M$ by facility $j \in N$. $K_j$ represents the capacity of facility $j\in N$ and each customer $i \in M$ has a demand $d_i$. Let $F_j$ be fixed cost for opening facility $j \in N$ and $\phi_j(z_j)$ be inverse S-shaped production cost of facility $j \in N$ for producing $z_j$ units. The deflection point/ economic point is denoted as $z_j^0$. We define decision variables $x_{ij}$ as the quantity transported from facility $j\in N$ to customer $i \in M$. Let $l_{j0}$ as binary variable that takes value one if facility $j\in N$ produces at or under economic point, and $l_{j1}$ as a binary variable that takes value one if the facility produces over economic point. If we define $\lambda_j$ as a binary variable that takes value 1 if facility $j\in N$ produces any positive quantity, then a facility location problem with inverse S-shaped production cost can be mathematically formulated as follows:
	\begin{align}
	\min_{\lambda,x,z} & \; \sum_{j=1}^{n}F_j \lambda_j +\sum_{i=1}^{m}\sum_{j=1}^{n}c_{ij}x_{ij} + \sum_{j=1}^{n} \phi_j(z_j) \label{eq:origObj6}\\
	\st & \sum_{j=1}^{n}x_{ij} \geq d_i && i=1,\dots,m \label{eq:origConstt1}\\
	& \lambda_j =l_{j0}+l_{j1}  && j=1,\dots,n  \label{eq:origConstt2}\\
	& \lambda_j \leq 1 && j=1,\dots,n  \label{eq:origConstt3}\\
	& z_j^0l_{j1}\leq z_j \leq K_jl_{j1}+z_j^0 l_{j0} && j=1,\dots,n \label{eq:origConstt4}  \\
	& z_j \geq \sum_{i=1}^{m}x_{ij} &&  j=1,\dots,n \label{eq:origConstt5}\\
	& l_{j0},l_{j1} \in \{0,1\} && j=1,\dots,n, \label{eq:origConstt6}\\
	& x_{ij} \geq 0 &&i=1,\dots,m, j=1;\dots,n, \label{eq:origConstt7}
	\end{align}
The above program is a special case of \eqref{eq:origObj}-\eqref{eq:origConst5}. Constraints \eqref{eq:origConstt1} and \eqref{eq:origConstt4} are demand and capacity constraints, respectively. Constraints \eqref{eq:origConstt3} indicate that any facility can either be operating below economic point $z_j^0$ or above it. $\phi_j(z_j)$ is a continuous inverse S-shaped function but not necessarily differentiable. 
\section{Computational Experiments} \label{sec:CompExper}
In this section, we discuss the data generation scheme, followed by the computational results. All computational experiments are conducted on a server with Intel(R) Xenon(R) Gold Dual-core processor CPU 6240 $@2.60$ GHz and 64 GB RAM. As explained in section~\ref{sec:algorithmDescription}, the IOA algorithm is coded in C++, and the MIP in Step 5 of the cutting plane algorithm (discussed in Appendix C) is solved using standard Branch$\&$Cut Solver (CPLEX 12.7.1). The optimality gap ($\epsilon$) is considered as $\frac{UB-LB}{LB} \times 100$, where UB and LB indicate the upper and lower bound of the original problem, respectively. The termination conditions for the algorithm are set to be $\epsilon =0.01$ in step 4 of algorithm~\ref{IOAalgo} or a maximum computational time of 3 hours, whichever reaches first. We aim to compare the computational  performance of our algorithm against the methods reported by \cite{lu2014facility}. The computational experiments by \cite{lu2014facility} were performed on a PC with Xeon 3.2 GHz and 8 GB RAM. Since \cite{lu2014facility} did not mention the core details of the processor, we report our computational results along with their reported computational results. 
\subsection{Data synthesis scheme} \label{subsec:datascheme}
The datasets are based on 71 facility location instances reported by \cite{holmberg1999exact}. \cite{holmberg1999exact} randomly generated four sets of test problems with different sizes and properties. \cite{lu2014facility} used three sets of the test problems out of four, i.e., 55 instances out of 71, and we consider those three sets in our computational experiments. The data characteristics are shown in Table 6 (see Appendix D), where $K/D$ indicates the ratio of the total capacity to the total demand. The coordinates of the locations are denoted by $L$, and the Euclidean distance between facility $j \in N$ and customer $i \in M$ is denoted by $e_{ij}$. The p-median problem \citep{klincewicz1990fleet} is used to select the locations of the potential facilities from the set of coordinates. The data synthesis strategy used for facility $j \in N$ and customer $i \in M$ by \cite{holmberg1999exact} is described briefly in Appendix D. 

The problem studied by \cite{holmberg1999exact} is a single sourcing facility location problem. \cite{lu2014facility} modified the dataset proposed by \cite{holmberg1999exact} slightly to accommodate it for multiple sourcing facility location problems. \cite{holmberg1999exact} considers the transportation cost between facility $j \in N$ and customer $i \in M$ as the cost for fulfilling the entire customer demand $d_i$ in the context of single sourcing. \cite{lu2014facility} divides the transportation cost $c_{ij}$ by the complete customer demand $d_i \;\forall i \in M, j \in N$. We use the modified dataset in our study.


We use three different types of S-shaped functions for each of the 55 problems to test the algorithm. The first is 
\begin{equation*}
	\phi_j(z_j) = 
	\begin{cases}
	\alpha_1 (z_j)^{\beta_1} \qquad \qquad \qquad \qquad \qquad \; \; \;  \; 0 \leq z_j \leq z_j^0\\
	\alpha_1 (z_j^0)^{\beta_1}+ \alpha_2(z_j-z_j^0)^{\beta_2}  \qquad \qquad z_j^0 \leq z_j \leq K \label{eq:fundef4}
	\end{cases}       
\end{equation*}
The second is
\begin{equation*}
	\phi^j(z_j) = 
	\begin{cases}
	\alpha_j (z_j)^{\beta_1} \qquad \qquad \qquad \qquad \; \;\;\;  \; \; \;  0 \leq z_j \leq z_j^0\\
	\alpha_1 (z_j^0)^{\beta_1}+ \alpha_2\frac{(z_j-z_j^0)}{(K_j-z_j)}   \qquad \qquad z_j^0 \leq z_j < K_j \label{eq:fundef5}
	\end{cases}       
\end{equation*}
These two function types are continuous but not necessarily differentiable at the economic point $x^0$. The third function is cubic function such as $ax^3-bx^2+cx_d$ where $a,b,c,d \geq 0$. Alternatively the function can be presented as below.
\begin{equation*}
    \phi_j(z_j)=a(z_j-z_j^0)^3 +\epsilon(z_j-z_j^0)+w, \;a>0, \;z_j^0>0, \;\epsilon \geq 0
\end{equation*}
Since the function $\phi^j(z_j)$ has to be monotonically increasing, $a,\epsilon$ must be chosen carefully so that $3a(z_j-z_j^0)^2 +\epsilon$ has at most one solution. Additionally, to test the efficiency of our algorithm, we use four different cost structures for all three function types leading to 12 test cases for each problem. \\
The cost structures are discussed below.
\begin{enumerate}[(i)]
    \item Base case: The three cost components (production cost, transportation cost, fixed cost) share the same percentage of the total cost.
    \item Dominant fixed cost case: The fixed cost dominates transportation cost and production cost.
    \item Dominant production cost case: The production cost dominates the rest two.
    \item Dominant transportation cost case: The transportation cost dominates the rest two.
\end{enumerate}
Table \ref{tab:funtype} depicts the cost setting for the three function types. We have considered $z_j^0=0.5K_j$. Fifty five data instances were solved for 12 function types resulting in a total of 660 test cases.
\begin{table}[htbp]
  \centering
  \caption{Cost Setting for different function types}
    \resizebox{\columnwidth}{!}{%
    \begin{tabular}{lrlrl}
    \toprule
    \multicolumn{1}{l}{Cost structure} &       & Production cost function &       & Other costs Setting \\
          &       &       &       &  \\
          \midrule
          \multicolumn{5}{c}{Function type 1} \\
                    &       &       &       &  \\
    1     &       & \begin{math}
    \phi_j(z_j)=
  \left\{
    \begin{array}{l}
	40 (z_j)^{0.5}  \qquad \qquad \qquad  \; \; \; \; \;\;\;\; 0 \leq z_j \leq z_j^0,\;j\in \{1,n\}\\
	40 (z_j^0)^{0.5}+ 0.5(z_j-z_j^0)^{2}   \;\;\;\; z_j^0 \leq z_j \leq K_j ,\;j\in \{1,n\}
\end{array}
  \right.
\end{math}     &       & Original  \\
    &         &$z_j^0= 0.5 \; K_j, \; j\in \{1,n\}$     &       & \\
            &       &     &       & \\
    2     &       & Same as 1     &       & Fixed cost is increased 10 fold for each facility. \\
                &       &     &       & \\
    3     &       & \begin{math}
    \phi_j(z_j)=
  \left\{
    \begin{array}{l}
	400 (z_j)^{0.5} \qquad \qquad \qquad   \; \;\;\;\; 0 \leq z_j \leq z_j^0,\;j\in \{1,n\}\\
	400 (z_j^0)^{0.5}+ 5(z_j-z_j^0)^{2}  \;\;\;\; z_j^0 \leq z_j \leq K_j ,\;j\in \{1,n\}
	\end{array}
  \right.
\end{math} &       & Original \\
                &        &$z_j^0= 0.5 \; K_j, \;  j\in \{1,n\}$     &       & \\
            &       &     &       & \\
    4     &       & Same as 1     &       & The Transportation cost is increased 10 fold between  \\
          &       &     &       &  source to destination across all locations. \\
                \midrule
          \multicolumn{5}{c}{Function type 2} \\
        &       &     &       &   \\
         1     &       & \begin{math}
    \phi_j(z_j)=
  \left\{
    \begin{array}{l}
	40 (z_j)^{0.5}  \qquad \qquad \qquad \qquad  \;\;\; \; \; \; \; \;\;\;\;  0 \leq z_j \leq z_j^0 \\
	40 (z_j^0)^{0.5}+ 20 (z_j^0)^{0.5}\frac{(z_j-z_j^0)}{(K_j-z_j)}   \qquad z_j^0 \leq z_j \leq K_j
\end{array}
  \right.
\end{math}     &       & Original  \\
    &         &$z_j^0= 0.5 \; K_j, \; j\in \{1,n\}$     &       & \\
            &       &     &       & \\
    2     &       & Same as 1     &       & Fixed cost is increased 10 fold for each facility. \\
                &       &     &       & \\
    3     &       & \begin{math}
    \phi_j(z_j)=
  \left\{
    \begin{array}{l}
	400 (z_j)^{0.5} \qquad \qquad \qquad \qquad  \qquad \;\;\;  0 \leq z_j \leq z_j^0,\\
	400 (z_j^0)^{0.5}+ 200 (z_j^0)^{0.5}\frac{(z_j-z_j^0)}{(K_j-z_j)}   \qquad z_j^0 \leq z_j \leq K_j 
	\end{array}
  \right.
\end{math} &       & Original  \\
                &        &$z_j^0= 0.5 \; K_j, \;  j\in \{1,n\}$     &       & \\
            &       &     &       & \\
    4     &       & Same as 1     &       & The transportation cost is increased 10 fold between  \\
          &       &     &       &  source to destination across all locations. \\
          \midrule
          \multicolumn{5}{c}{Function type 3} \\
        &       &     &       &   \\
            1     &       & $\phi(z_j)=10^{-4}(z_j)^3 -3 \times 10^{-4}z_j^0(z_j)^2+3 \times 10^{-4}(z_j^0)^2z_j$     &       & Original  \\
                 &         &$z_j^0= 0.5 \; K_j, \; j\in \{1,n\}$     &       & \\
            &       &     &       & \\
    2     &       & Same as 1     &       & Fixed cost is increased 10 fold for each facility. \\
                &       &     &       & \\
    3   &       & $\phi_j(z_j)=10^{-3}(z_j)^3 -3 \times 10^{-3}z_j^0(z_j)^2+3 \times 10^{-3}(z_j^0)^2 z_j$     &       & Original  \\
                &        &$z_j^0= 0.5 \; K_j, \;  j\in \{1,n\}$     &       & \\
            &       &     &       & \\
    4     &       & Same as 1     &       & The transportation cost is increased 10 fold between  \\
          &       &     &       &  source to destination across all locations. \\
          \bottomrule
    \end{tabular}%
  \label{tab:funtype}%
  }
\end{table}%
\subsection{Computational Results}	\label{sec:results1}
In this section, we compare the computational performance of the proposed algorithm against the method reported by \cite{lu2014facility}. \cite{lu2014facility} proposed a method based on Lagrangian relaxation, column generation (CG), and branch-and-bound. \cite{lu2014facility} proposed three CG heuristics based on the extent of branching. They are CG heuristic (CG0), Enhanced CG heuristic (CG1), and Branch-and-price heuristic (CG2). We compare our proposed IOA approach against the three heuristics. The proposed IOA method is exact in nature and hence guarantees global optimality on termination. The method converges to global optimality for 95 percent of the 660 test cases within the maximum allowed CPU Time of 3 hours. For the rest of the test instances, the IOA algorithm terminates at an average optimality gap of 0.55$\%$. Since \cite{lu2014facility} reported only the average, minimum, and maximum of both computational time and the optimality gap ($\%$) across various problems for each cost function and cost structure type, we also do the same for a meaningful comparison. The results are reported through Tables~\ref{tab:basecase}-\ref{tab:domvarcase}, where the best optimality gap ($\%$) among the four methods is highlighted in bold for every cost function type and cost structure type, while the best CPU time among the four methods is highlighted in bold only if its corresponding gap is also the best among the four. Tables~\ref{tab:basecase}-\ref{tab:domvarcase} show the performance comparisons of the proposed IOA method against CG0, CG1, and CG2, corresponding to the three different functions and the four different cost structures. The missing values in the rows indicate that the authors did not report the corresponding results. Tables~\ref{tab:basecase}-\ref{tab:domvarcase} clearly show that the IOA algorithm outperforms CG0, CG1, and CG2 for most of the test cases. In fact, for 60 percent of the test cases, the IOA algorithm outperforms all CG0, CG1, and CG2 while reporting a 0.00\% optimality gap. Next, we describe the other observations briefly from the tables. \par

It is evident from the results that the IOA algorithm is able to converge to the global optimum with an average gap of 0.00\% in all the test datasets when $n=10$ for all function types and cost structures, whereas CG0 and CG2 methods are unable to converge to the global optimum in all these test cases. The results for CG1 and CG2 are not reported by the author for most of the test cases. It is important to note that a maximum optimality gap (maximum among all functions and cost structure types when $n=10$) of 11.52 percent is observed in both CG0 and CG2 for the function type 3 for the dominant fixed cost structure. Similarly, for the dominant transportation cost case, the IOA algorithm outperforms on average all cases by converging to the global optimum with lower computational time than both CG0 and CG2 for all function types when $n=10$. Though the IOA algorithm takes a higher average CPU time than both CG0 and CG2 for some test cases, both CG0 and CG2 report a higher average optimality gap $(\%)$ than the IOA algorithm.

For instances with $n \leq 20$ and $n \leq 30$ as well, the IOA algorithm converges to the global optimum for most of the data instances corresponding to all function types and cost structures. The IOA algorithm reports an average gap of either 0 or very close to 0 for all function types and cost structures except the dominant production cost case of function 1. For the dominant production cost case of function 1, the average optimality gaps for datasets with $n\leq20$ and $n\leq30$ at termination are found to be 0.28 percent and 0.51 percent, respectively. 

\begin{table}[htbp]
\centering
    \caption{Experimental results for base case}
       \begin{threeparttable}
    \begin{tabular}{ccccccccccccccc}
    \toprule
   &     &     &     &     & \multicolumn{4}{c}{Optimality Gap (\%)} &     &     & \multicolumn{4}{c}{CPU Time (seconds)} \\
\cmidrule{5-8}\cmidrule{11-14}    n   &    m   &     &     & IOA &  CG0 & CG1 & CG2 &     &     & IOA &  CG0 & CG1 & CG2 \\
    \midrule
       &    & \multicolumn{12}{c}{Function Type 1} \\
    =10 & =50 & Avg &    & \textbf{0.00} & 0.68 & -  & 0.68 &    &    & 36.97 & 1.34 & -  & 1.96 \\
       &    & Min &    & \textbf{0.00} & 0.00 & -  & 0.00 &    &    & \textbf{0.14} & 0.98 & -  & 1.16 \\
       &    & Max &    & \textbf{0.00} & 3.07 & -  & 3.07 &    &    & 154.18 & 2.20 & -  & 3.85 \\
    =10 & $\leq$100 & Avg &    & 0.00 & -  & -  & -  &    &    & 25.46 & - & -  & - \\
       &    & Min &    & 0.00 & -  & -  & -  &    &    & 0.14 & -  & -  & - \\
       &    & Max &    & 0.00 & -  & -  & -  &    &    & 154.18 & - & -  & - \\
      &    &    &    &    &    &    &    &    &    &    &    &    &  \\
    $\leq$20 & $\leq$100 & Avg &    & \textbf{0.03} & - & 0.66 & -  &    &  & 1409.73 & -  & 516.73 & - \\
       &    & Min &    & \textbf{0.00} & - & 0.00 & -  &    &  & \textbf{0.14} & -  & 7.22 & - \\
       &    & Max &    & \textbf{0.52} & -  & 2.65 & -  &    &  & 10800.00 & -  & 4988.04 & - \\
       &    &    &    &    &    &    &    &    &  &    &    &    &  \\
    $\leq$30 & $\leq$150 & Avg &    & \textbf{0.08} & 0.86 & -  & -  &    &  & 2857.41 & 670 & -  & - \\
       &    & Min &    & \textbf{0.00} & 0.00 & -  & -  &    &  & \textbf{0.14} & 2.70 & -  & - \\
       &    & Max &    & \textbf{1.42} & 2.65 & -  & -  &    &  & 10800.00 & 3205.79 & -  & - \\
    \midrule
       &    & \multicolumn{12}{c}{Function Type 2} \\
    =10 & =50 & Avg &    & \textbf{0.00} & 0.97 & -  & 0.92 &    &    & 3.64 & 2.20 & -  & 10.69 \\
       &    & Min &    & \textbf{0.00} & 0.00 & -  & 0.00 &    &    & \textbf{0.29} & 0.98 & -  & 1.16 \\
       &    & Max &    & \textbf{0.00} & 2.71 & -  & 2.29 &    &    & 14.37 & 5.07 &    & 37.98 \\
    =10 & $\leq$100 & Avg &    & 0.00 & -  & -  & -  &    &    & 3.28 & -  & -  & - \\
       &    & Min &    & 0.00 & -  & -  & -  &    &    & 0.29 & -  & -  & - \\
       &    & Max &    & 0.00 & -  & -  & -  &    &    & 14.37 & -  & -  & - \\
  &    &    &    &    &    &    &    &    &    &    &    &    &  \\
    $\leq$20 & $\leq$100 & Avg &    & \textbf{0.00} & - & 0.81 & -  &    &  & \textbf{103.48} & - & 449.54 & - \\
       &    & Min &    & \textbf{0.00} & - & 0.00 & -  &    &  & \textbf{0.29} & - & 17.88 & - \\
       &    & Max &    & \textbf{0.00} & - & 2.95 & -  &    &  & \textbf{1362.60} & - & 2587.74 & - \\
       &    &    &    &    &    &    &    &    &  &    &    &    &  \\
    $\leq$30 & $\leq$150 & Avg &    & \textbf{0.00} & 0.93 & -  & -  &    &  & \textbf{95.89} & 1485.21 & -  & - \\
       &    & Min &    & \textbf{0.00} & 0.00  & -  & -  &    &  & \textbf{0.29} & 3.74 & -  & - \\
       &    & Max &    & \textbf{0.00} & 2.95 & -  & -  &    &  & \textbf{1362.60} & 13707.28 & -  & - \\
    \midrule
       &    & \multicolumn{12}{c}{Function Type 3} \\
    =10 & =50 & Avg &    & \textbf{0.00} & 0.62 & -  & 0.62 &    &    & \textbf{0.77} & 1.63 & -  & 5.28 \\
       &    & Min &    & \textbf{0.00} & 0.00  & -  & 0.00  &    &    & \textbf{0.15} & 0.98 & -  & 1.62 \\
       &    & Max &    & \textbf{0.00} & 2.72 & -  & 2.72 &    &    & 4.40 & 2.97 & -  & 10.25 \\
    =10 & $\leq$100 & Avg &    & 0.00 & -  & -  & -  &    &    & 0.86 & -  & -  & - \\
       &    & Min &    & 0.00 & -  & -  & -  &    &    & 0.15 & -  & -  & - \\
       &    & Max &    & 0.00 & -  & -  & -  &    &    & 4.40 & -  & -  & - \\
          &    &    &    &    &    &    &    &    &    &    &    &    &  \\
    $\leq$20 & $\leq$100 & Avg &    & \textbf{0.00} & - & 1.02 & -  &    &  & \textbf{49.03} & - & 283.00 & - \\
       &    & Min &    & \textbf{0.00} & - & 0.00 & -  &    &  & \textbf{0.15} & - & 4.46 & - \\
       &    & Max &    & \textbf{0.00} & - & 2.94 & -  &    &    & \textbf{530.78} & - & 3098.23 & - \\
       &    &    &    &    &    &    &    &    &    &    &    &    &  \\
    $\leq$30 & $\leq$150 & Avg &    & \textbf{0.00} & 1.23 & -  & -  &    &    & \textbf{461.33} & 1796.57 & -  & - \\
       &    & Min &    & \textbf{0.00} & 0.00 & -  & -  &    &    & \textbf{0.15} & 3.49 & -  & - \\
       &    & Max &    & \textbf{0.00} & 5.03 & -  & -  &    &    & \textbf{9493.66} & 16680.61 & -  & - \\
\bottomrule  
\end{tabular}%
\begin{tablenotes}[normal,flushleft]
\item  $-$ indicates that the result is not reported by \cite{lu2014facility} for the respective method.
\end{tablenotes}
 \end{threeparttable}
\label{tab:basecase}%
\end{table}%

\begin{table}[htbp]
  \centering
  \caption{Experimental results for dominant fixed cost case}
    \begin{threeparttable}
    \begin{tabular}{ccccccccccccccc}
    \toprule
      &     &     &     &     & \multicolumn{4}{c}{Optimality Gap (\%)} &     &     & \multicolumn{4}{c}{CPU Time (seconds)} \\
\cmidrule{5-8}\cmidrule{11-14}    n   &    m   &     &     & IOA &  CG0 & CG1 & CG2 &     &     & IOA &  CG0 & CG1 & CG2 \\
    \midrule
       &    & \multicolumn{11}{c}{Function Type 1}                 &  \\
    =10 & =50 & Avg &    & \textbf{0.00} & 0.47 & -  & 0.46 &    &    & \textbf{0.79} & 1.72 & -  & 14.44 \\
       &    & Min &    & \textbf{0.00} & 0.00 & -  & 0.00 &    &    & \textbf{0.13} & 1.00 & -  & 1.15 \\
       &    & Max &    & \textbf{0.00} & 1.73 & -  & 1.72 &    &    & 4.98 & 2.92 & -  & 115.86 \\
    =10 & $\leq$100 & Avg &    & 0.00 & -  & -  & -  &    &    & 0.70 & -  & -  & - \\
       &    & Min &    & 0.00 & -  & -  & -  &    &    & 0.13 & -  & -  & - \\
       &    & Max &    & 0.00 & -  & -  & -  &    &    & 4.98 & -  & -  & - \\
       &    &    &    &    &    &    &    &    &  &    &    &    &  \\
    $\leq$20 & $\leq$100 & Avg &    & \textbf{0.00} & - & 0.59 & -  &    &  & \textbf{32.11} & - & 1151.18 & - \\
       &    & Min &    & \textbf{0.00} & - & 0.00 & -  &    &  & \textbf{0.13} & - & 10.36 & - \\
       &    & Max &    & \textbf{0.00} & - & 2.49 & -  &    &  & \textbf{678.86} & - & 22098.06 & - \\
       &    &    &    &    &    &    &    &    &  &    &    &    &  \\
    $\leq$30 & $\leq$150 & Avg &    & \textbf{0.00} & 0.78 & -  & -  &    &  & \textbf{194.21} & 1177.83 & -  & - \\
       &    & Min &    & \textbf{0.00} & 0.00 & -  & -  &    &  & \textbf{0.13} & 3.09 & -  & - \\
       &    & Max &    & \textbf{0.00} & 3.98 & -  & -  &    &  & \textbf{3018.88} & 16872.55 & -  & - \\
    \midrule
       &    & \multicolumn{11}{c}{Function Type 2}                 &  \\
    =10 & =50 & Avg &    & \textbf{0.00} & 3.68 & -  & 2.89 &    &    & \textbf{0.47} & 3.12 & -  & 485.00 \\
       &    & Min &    & \textbf{0.00} & 0.28 & -  & 0.28 &    &    & \textbf{0.26} & 1.00 & -  & 3.83 \\
       &    & Max &    & \textbf{0.00} & 7.24 & -  & 6.63 &    &    & \textbf{1.12} & 11.58 & -  & 3486.1 \\
    =10 & $\leq$100 & Avg &    & 0.00 & -  & -  & -  &    &    & 0.45 & -  & -  & - \\
       &    & Min &    & 0.00 & -  & -  & -  &    &    & 0.26 & -  & -  & - \\
       &    & Max &    & 0.00 & -  & -  & -  &    &    & 1.12 & -  & -  & - \\
       &    &    &    &    &    &    &    &    &  &    &    &    &  \\
    $\leq$20 & $\leq$100 & Avg &    & \textbf{0.00} & - & 2.19 & -  &    &  & \textbf{2.96} & - & 1601.85 & - \\
       &    & Min &    & \textbf{0.00} & - & 0.13 & -  &    &  & \textbf{0.26} & - & 10.05 & - \\
       &    & Max &    & \textbf{0.00} & - & 6.26 & -  &    &  & \textbf{43.04} & - & 29853.99 & - \\
       &    &    &    &    &    &    &    &    &  &    &    &    &  \\
    $\leq$30 & $\leq$150 & Avg &    & \textbf{0.00} & 2.67 & -  & -  &    &  & \textbf{7.73} & 2153.22 & -  & - \\
       &    & Min &    & \textbf{0.00} & 0.00 & -  & -  &    &  & \textbf{0.26} & 2.95 & -  & - \\
       &    & Max &    & \textbf{0.00} & 11.67 & -  & -  &    &  & \textbf{132.63} & 36122.19 & -  & - \\
    \midrule
       &    & \multicolumn{11}{c}{Function Type 3}                 &  \\
    =10 & =50 & Avg &    & \textbf{0.00} & 5.30 & -  & 4.36 &    &    & \textbf{0.35} & 1.37 & -  & 61.65 \\
       &    & Min &    & \textbf{0.00} & 0.00 & -  & 0.00 &    &    & \textbf{0.09} & 0.98 & -  & 7.74 \\
       &    & Max &    & \textbf{0.00} & 11.52 & -  & 11.52 &    &    & \textbf{0.71} & 2.94 & -  & 271.45 \\
    =10 & $\leq$100 & Avg &    & 0.00 & -  & -  & -  &    &    & 0.34 & -  & -  & - \\
       &    & Min &    & 0.00 & -  & -  & -  &    &    & 0.09 & -  & -  & - \\
       &    & Max &    & 0.00 & -  & -  & -  &    &    & 0.71 & -  & -  & - \\
       &    &    &    &    &    &    &    &    &  &    &    &    &  \\
    $\leq$20 & $\leq$100 & Avg &    & \textbf{0.00} & - & 3.4 & -  &    &  & \textbf{2.68} & - & 464.54 & - \\
       &    & Min &    & \textbf{0.00} & - & 0.00 & -  &    &  & \textbf{0.09} & - & 2.11 & - \\
       &    & Max &    & \textbf{0.00} & - & 10.04 & -  &    &  & \textbf{25.25} & - & 2593.05 & - \\
       &    &    &    &    &    &    &    &    &  &    &    &    &  \\
    $\leq$30 & $\leq$150 & Avg &    & \textbf{0.00} & 4.08 & -  & -  &    &  & \textbf{14.20} & 1303.75 & -  & - \\
       &    & Min &    & \textbf{0.00} & 0.00 & -  & -  &    &  & \textbf{0.09} & 2.20 & -  & - \\
       &    & Max &    & \textbf{0.00} & 21.20 & -  & -  &    &    & \textbf{213.32} & 15972.76 & -  & - \\
\bottomrule
\end{tabular}%
\begin{tablenotes}[normal,flushleft]
\item   $-$ indicates that the result is not reported by \cite{lu2014facility} for the respective method.
\end{tablenotes}
  \label{tab:fixedcostcase}%
 \end{threeparttable}
\end{table}%

\begin{table}[htbp]
  \centering
  \caption{Experimental results for dominant production cost case}
    \begin{tabular}{ccccccccccccccc}
    \toprule
      &     &     &     &     & \multicolumn{4}{c}{Optimality Gap (\%)} &     &     & \multicolumn{4}{c}{CPU Time (seconds)} \\
\cmidrule{5-8}\cmidrule{11-14}    n   &    m   &     &     & IOA &  CG0 & CG1 & CG2 &     &     & IOA &  CG0 & CG1 & CG2 \\
\midrule
       &    & \multicolumn{12}{c}{Function type 1} \\
    =10 & =50 & Avg &    & \textbf{0.00} & 0.08 & -  & 0.08 &    &    & 1325.10 & 1.33 & -  & 3.30 \\
       &    & Min &    & \textbf{0.00} & 0.00 & -  & 0.00 &    &    & \textbf{0.14} & 0.98 & -  & 1.17 \\
       &    & Max &    & \textbf{0.00} & 0.36 & -  & 0.36 &    &    & 5643.72 & 1.84 & -  & 5.34 \\
    \multicolumn{1}{c}{=10} & \multicolumn{1}{c}{$\leq$100} & Avg &    & 0.00 & -  & -  & -  &    &    & 884.69 & -  & -  & - \\
       &    & Min &    & 0.00 & -  & -  & -  &    &    & 0.14 & -  & -  & - \\
       &    & Max &    & 0.00 & -  & -  & -  &    &    & 5643.72 & -  & -  & - \\
       &    &    &    &    &    &    &    &    &    &    &    &    &  \\
    $\leq$20 & $\leq$100 & Avg &    & \textbf{0.28} & - & 0.42 & -  &    &  & 2801.59 & -  & 1339.96 & - \\
       &    & Min &    & \textbf{0.00} & - & 0.00 & -  &    &  & \textbf{0.14} & -  & 10.33 & - \\
       &    & Max &    & 2.29 & - & \textbf{1.33} & -  &    &  & 10800.00 & -  & 17864.00 & - \\
       &    &    &    &    &    &    &    &    &  &    &    &    &  \\
    $\leq$30 & $\leq$150 & Avg &    & 0.51 & \textbf{0.48} & -  & -  &    &  & 5564.68 & \textbf{984.70}  & - & - \\
       &    & Min &    & \textbf{0.00} & 0.00 & -  & -  &    &  & \textbf{0.14} & 3.63  & - & - \\
       &    & Max &    & 2.29 & \textbf{2.26} & -  & -  &    &  & 10800.00 & \textbf{5085.60}  & - & - \\
    \midrule
       &    & \multicolumn{12}{c}{Function type 2} \\
    =10 & =50 & Avg &    & \textbf{0.00} & 0.50 & -  & 0.50 &    &    & 4.34 & 1.67 & -  & 24.94 \\
       &    & Min &    & \textbf{0.00} & 0.00 & -  & 0.00  &    &    & \textbf{0.49} & 1.08 & -  & 1.65 \\
       &    & Max &    & \textbf{0.00} & 1.71 & -  & 1.70 &    &    & 9.93 & 2.38 & -  & 101.43 \\
    \multicolumn{1}{c}{=10} & \multicolumn{1}{c}{$\leq$100} & Avg &    & 0.00 & -  & -  & -  &    &    & 4.37 & -  & -  & - \\
       &    & Min &    & 0.00 & -  & -  & -  &    &    & 0.49 & -  & -  & - \\
       &    & Max &    & 0.00 & -  & -  & -  &    &    & 9.93 & -  & -  & - \\
       &    &    &    &    &    &    &    &    &  &    &    &    &  \\
    $\leq$20 & $\leq$100 & Avg &    & \textbf{0.00} & - & 0.38 & -  &    &  & \textbf{136.30} & -  & 409.11 & - \\
       &    & Min &    & \textbf{0.00} & - & 0.00 & -  &    &  & \textbf{0.49} & -  & 13.74 & - \\
       &    & Max &    & \textbf{0.00} & - & 1.56 & -  &    &  & \textbf{876.43} & -  & 1839.91 & - \\
       &    &    &    &    &    &    &    &    &  &    &    &    &  \\
    $\leq$30 & $\leq$150 & Avg &    & \textbf{0.02} & 0.53 & -  & -  &    &  & 1531.95 & 810.04 & -  & - \\
       &    & Min &    & \textbf{0.00} & 0.00 & -  & -  &    &  & \textbf{0.49} & 4.49 & -  & - \\
       &    & Max &    & \textbf{1.34} & 2.72 & -  & -  &    &  & 10800.00 & 7152.37 & -  & - \\
    \midrule
       &    & \multicolumn{12}{c}{Function type 3} \\
    =10 & =50 & Avg &    & \textbf{0.00} & 0.86 & -  & 0.86 &    &    & 11.03 & 2.51 & -  & 7.36 \\
       &    & Min &    & \textbf{0.00} & 0.00 & -  & 0.00  &    &    & \textbf{0.58} & 1.05 & -  & 2.17 \\
       &    & Max &    & \textbf{0.00} & 2.06 & -  & 2.06 &    &    & 73.71 & 6.63 & -  & 26.77 \\
    =10 & $\leq$100 & Avg &    & 0.00 & -  & -  & -  &    &    & 7.65 & -  & -  & - \\
       &    & Min &    & 0.00 & -  & -  & -  &    &    & 0.51 & -  & -  & - \\
       &    & Max &    & 0.00 & -  & -  & -  &    &    & 73.71 & -  & -  & - \\
       &    &    &    &    &    &    &    &    &  &    &    &    &  \\
    $\leq$20 & $\leq$100 & Avg &    & \textbf{0.05} & - & 0.92 & -  &    &  & 1387.47 &    & 172.81 & - \\
       &    & Min &    & \textbf{0.00} & - & 0.00  & -  &    &  & \textbf{0.51} &    & 13.21 & - \\
       &    & Max &    & \textbf{1.02} & - & 5.97 & -  &    &  & 10800 &    & 859.50 & - \\
       &    &    &    &    &    &    &    &    &    &    &    &    &  \\
    $\leq$30 & $\leq$150 & Avg &    & \textbf{0.05} & 1.18 & -  & -  &    &    & 1828.86 & 1823.92 & -  & - \\
       &    & Min &    & \textbf{0.00} & 0.00  & -  & -  &    &    & \textbf{0.51} & 5.58 & -  & - \\
       &    & Max &    & \textbf{1.02} & 6.17 & -  & -  &    &    & \textbf{10800.00} & 19674.19 & -  & - \\
    \bottomrule   \end{tabular}%
       \begin{tablenotes}[normal,flushleft]
\item  $-$ indicates that the result is not reported by \cite{lu2014facility} for the respective method.
\end{tablenotes}
  \label{tab:domprodcase}%
\end{table}

\begin{figure}
        \centering
		\includegraphics[width=1\linewidth]{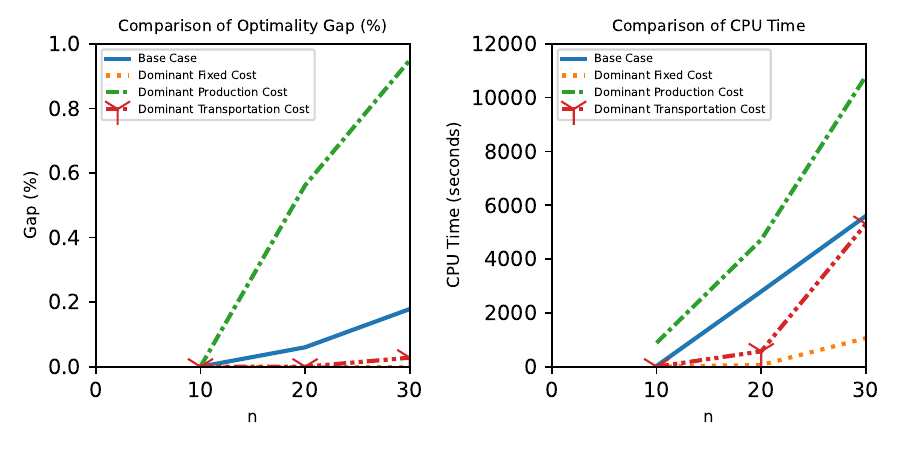}
		\caption{Computational performance of the IOA algorithm for function type 1 (n=10,20,30)}
		\label{fig:fun1}
	\end{figure}

\begin{figure}
		\centering
			\includegraphics[width=1\linewidth]{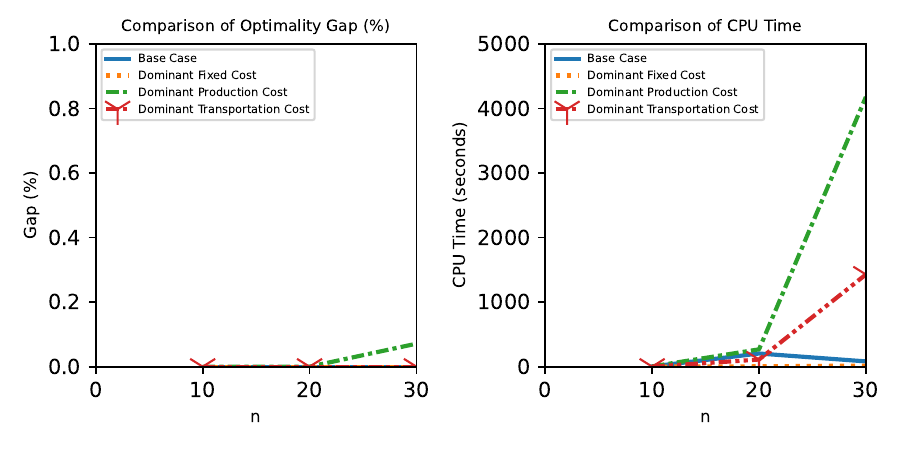}
		\caption{Computational performance of the IOA algorithm for function type 2 (n=10,20,30)}
		\label{fig:fun2}
	\end{figure}

\begin{figure}
		\centering
			\includegraphics[width=1\linewidth]{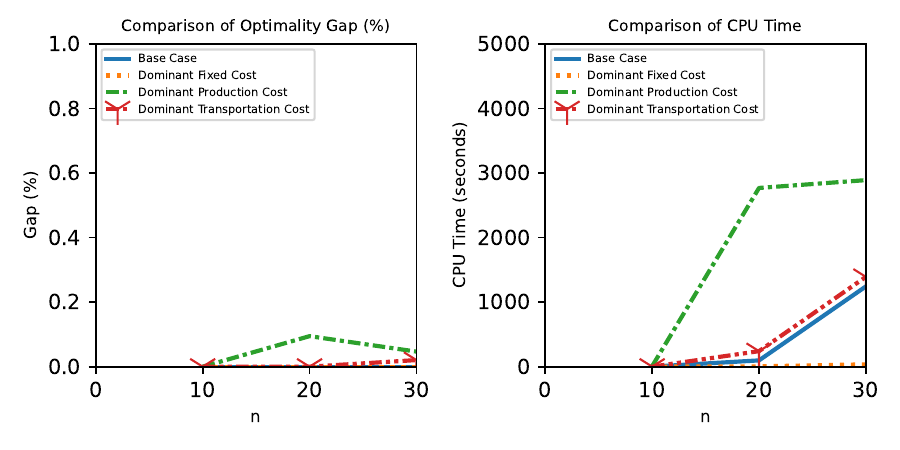}
		\caption{Computational performance of the IOA algorithm for function type 3 (n=10,20,30)}
		\label{fig:fun3}
	\end{figure}

\begin{table}[htbp]
  \centering
  \caption{Experimental results for dominant transportation cost case}
    \begin{tabular}{ccccccccccccccc}
    \toprule
      &     &     &     &     & \multicolumn{4}{c}{Optimality Gap (\%)} &     &     & \multicolumn{4}{c}{CPU Time (seconds)} \\
\cmidrule{5-8}\cmidrule{11-14}    n   &    m   &     &     & IOA &  CG0 & CG1 & CG2 &     &     & IOA &  CG0 & CG1 & CG2 \\
\midrule
       &    & \multicolumn{12}{c}{Function type 1} \\
    =10 & =50 & Avg &    & \textbf{0.00} & 0.00 & -  & 0.00 &    &    & \textbf{0.61} & 1.06 & -  & 1.09 \\
       &    & Min &    & \textbf{0.00} & 0.00 & -  & 0.00 &    &    & \textbf{0.16} & 0.98 & -  & 1.06 \\
       &    & Max &    & \textbf{0.00} & 0.01 & -  & 0.01 &    &    & 1.48 & 1.16 & -  & 1.16 \\
    \multicolumn{1}{c}{=10} & \multicolumn{1}{c}{$\leq$100} & Avg &    & 0.00 & -  & -  & -  &    &    & 0.72 & -  & -  & - \\
       &    & Min &    & 0.00 & -  & -  & -  &    &    & 0.16 & -  & -  & - \\
       &    & Max &    & 0.00 & -  & -  & -  &    &    & 2.31 & -  & -  & - \\
       &    &    &    &    &    &    &    &    &  &    &    &    &  \\
    $\leq$20 & $\leq$100 & Avg &    & \textbf{0.00} & - & 0.09 & -  &    &  & 284.71 & -  & 138.90 & - \\
       &    & Min &    & \textbf{0.00} & - & 0.00 & -  &    &  & \textbf{0.16} & -  & 8.31 & - \\
       &    & Max &    & \textbf{0.00} & - & 0.64 & -  &    &  & 1679.20 & -  & 1143.46 & - \\
       &    &    &    &    &    &    &    &    &  &    &    &    &  \\
    $\leq$30 & $\leq$150 & Avg &    & \textbf{0.01} & 0.23 & -  & -  &    &  & 2018.54 & 807.01 & -  & - \\
       &    & Min &    & \textbf{0.00} & 0.00 & -  & -  &    &  & \textbf{0.16} & 2.95 & -  & - \\
       &    & Max &    & \textbf{0.28} & 1.07 & -  & -  &    &  & 10800.00 & 6552.84 & -  & - \\
    \midrule
       &    & \multicolumn{12}{c}{Function type 2} \\
    =10 & =50 & Avg &    & \textbf{0.00} & 0.00 & -  & 0.00 &    &    & \textbf{0.48} & 1.12 & -  & 1.06 \\
       &    & Min &    & \textbf{0.00} & 0.00 & -  & 0.00 &    &    & \textbf{0.18} & 0.98 & -  & 1.01 \\
       &    & Max &    & \textbf{0.00} & 0.01 & -  & 0.00 &    &    & \textbf{0.77} & 1.31 & -  & 1.09 \\
    \multicolumn{1}{c}{=10} & \multicolumn{1}{c}{$\leq$100} & Avg &    & 0.00 & -  & -  & -  &    &    & 0.54 & -  & -  & - \\
       &    & Min &    & 0.00 & -  & -  & -  &    &    & 0.18 & -  & -  & - \\
       &    & Max &    & 0.00 & -  & -  & -  &    &    & 1.24 & -  & -  & - \\
       &    &    &    &    &    &    &    &    &  &    &    &    &  \\
    $\leq$20 & $\leq$100 & Avg &    & \textbf{0.00} & - & 0.06 & -  &    &  & \textbf{55.37} & -  & 181.77 & - \\
       &    & Min &    & \textbf{0.00} & - & 0.00 & -  &    &  & \textbf{0.18} & -  & 15.21 & - \\
       &    & Max &    & \textbf{0.00} & - & 0.45 & -  &    &  & \textbf{784.10} & -  & 992.71 & - \\
       &    &    &    &    &    &    &    &    &  &    &    &    &  \\
    $\leq$30 & $\leq$150 & Avg &    & \textbf{0.00} & 0.15 & -  & -  &    &  & 2018.54 & 791.68 & -  & - \\
       &    & Min &    & \textbf{0.00} & 0.00 & -  & -  &    &  & \textbf{0.16} & 5.65 & -  & - \\
       &    & Max &    & \textbf{0.00} & 1.30 & -  & -  &    &  & 10800.00 & 7021.46 & -  & - \\
    \midrule
       &    & \multicolumn{12}{c}{Function type 3} \\
    =10 & =50 & Avg &    & \textbf{0.00} & 0.00 & -  & 0.00 &    &    & \textbf{0.31} & 1.22 & -  & 1.10 \\
       &    & Min &    & \textbf{0.00} & 0.00 & -  & 0.00 &    &    & \textbf{0.18} & 1.02 & -  & 1.00 \\
       &    & Max &    & \textbf{0.00} & 0.00 & -  & 0.00 &    &    & \textbf{0.63} & 1.71 & -  & 1.43 \\
    \multicolumn{1}{c}{=10} & \multicolumn{1}{c}{$\leq$100} & Avg &    & 0.00 & -  & -  & -  &    &    & 0.34 & -  & -  & - \\
       &    & Min &    & 0.00 & -  & -  & -  &    &    & 0.09 & -  & -  & - \\
       &    & Max &    & 0.00 & -  & -  & -  &    &    & 0.71 & -  & -  & - \\
       &    &    &    &    &    &    &    &    &  &    &    &    &  \\
    $\leq$20 & $\leq$100 & Avg &    & \textbf{0.00} & - & 0.10 & -  &    &  & \textbf{120.24} & -  & 135.29 & - \\
       &    & Min &    & \textbf{0.00} & - & 0.00 & -  &    &  & \textbf{0.14} & -  & 9.81 & - \\
       &    & Max &    & \textbf{0.00} & - & 0.79 & -  &    &  & 1693.04 & -  & 661.55 & - \\
       &    &    &    &    &    &    &    &    &    &    &    &    &  \\
    $\leq$30 & $\leq$150 & Avg &    & \textbf{0.01} & 0.21 & -  & -  &    &    & \textbf{559.103} & 717.54 & -  & - \\
       &    & Min &    & \textbf{0.00} & 0.00 & -  & -  &    &    & \textbf{0.14} & 4.35 & -  & - \\
       &    & Max &    & \textbf{0.39} & 1.68 & -  & -  &    &    & 10800.000 & 5569.64 & -  & - \\
   \bottomrule   \end{tabular}%
       \begin{tablenotes}[normal,flushleft]
\item $-$ indicates that the result is not reported by \cite{lu2014facility} for the respective method.
\end{tablenotes}
  \label{tab:domvarcase}%
\end{table}
Next, we study the performance of IOA with the increase in the size of the problem. We plot the average computational time and average optimality gap ($\%$) against $n$ for each function type and cost structure type in Figures \ref{fig:fun1}-\ref{fig:fun3}. As shown in Figures \ref{fig:fun1}-\ref{fig:fun3}, with the increasing value of n, the average CPU time increases, but the rise of CPU time is not significant except in the dominant production cost case. 

\subsection{Analysis of Impact of Different Cost Components and Economic Point}	\label{sec:managementinsghts}
In this section, we examine the impact of different cost components and economic points for function type 2 on the design of the network and find some important managerial insights. Specifically, we change the weights for production cost, transportation cost, and the economic point and study how the solution changes. The results are presented through Tables 7-9 in Appendix E, where columns indicated as $n_e,\; n_d, and \; n_T$ indicate the total number of facilities operating under economies of scale, total number of facilities operating under dis-economies of scale, and total number of facilities open, respectively. The column, Obj/Gap $(\%)$, indicates the optimal objective function value when we are able to close the optimality gap $(\%)$ within the set computational time; otherwise, it indicates the optimality gap $(\%)$.

\subsubsection{Impact of Transportation Cost Between Facilities and Customers}	\label{sec:transportimpact}
We examine the impact of transportation costs on the design of the network. We define a weight $\alpha$ as the weight of the transportation cost and describe the total cost to meet the demand of the customers as follows,
$$
\mbox{Total cost}=\sum_{j=1}^{n}F_j \lambda_j +\alpha \sum_{i=1}^{m}\sum_{j=1}^{n} c_{ij}x_{ij} + \sum_{j=1}^{n} \phi_j(z_j)
$$
This total cost is the same as the base case when $\alpha=1$, whereas it represents the dominant transportation cost case when $\alpha=10$. To analyze the effect of $\alpha$ on the solution of the model, we vary $\alpha$ from 1 to 20 and report the results in Table 7 (see Appendix E), where we observe the following.

\paragraph{Observation 1:} The number of open facilities tends to increase with the increase in the weight $(\alpha)$\par
From Table 7, it is observable that the number of open facilities increases substantially across all dimensions when we increase the weight from 1 to 10. This behavior can be explained by the fact that customers tend to fulfill their demands from the nearest facilities when transportation costs are increased significantly. An increase in weight beyond 10 does not lead to a significant rise in open facilities as a high transportation cost has already forced the customers to be assigned to the nearest facilities.

\paragraph{Observation 2:} The number of open facilities operating under economies of scale tends to increase with the increase of the weight $(\alpha)$, whereas the number of open facilities operating under dis-economies of scale tends to decrease.\par
As shown in Table 7, the number of open facilities operating under economies of scale increases when the transportation costs are increased ten-fold $(\alpha=10)$. The reason for such a behavior is that for the base case $(\alpha=1)$, the model always tries to select facilities with large capacity to benefit from economies of scale, but when we increase the transportation cost, it outweighs the cost saving happening due to economies of scale. As a result, the model tries to assign customers to the nearest facilities. A further increase in $\alpha$ from 10 to 20, does not lead to a significant change for the same reason as explained in observation 1.

\subsubsection{Impact of Production Cost of the Facilities}	\label{sec:prodimpact}
Next, we examine the impact of the production costs on the network. Here, we define a weight $\theta$ for the production cost as follows,
$$
\mbox{Total cost} =\sum_{j=1}^{n}F_j \lambda_j + \sum_{i=1}^{m}\sum_{j=1}^{n} c_{ij}x_{ij} + \theta \sum_{j=1}^{n} \phi_j(z_j)
$$
To analyze how the values of $\theta$ affect the solution of the model, we vary $\theta$ from 1 to 20 and report the results in Table 8 (see Appendix E).\par
From Table 8, it is observable that the change in the number of open facilities is not significant when we increase $\theta$ from 1 to 20 for all data instances. As mentioned earlier, when $\theta=1$, the model tries to allocate customers to facilities with higher capacities to benefit from economies of scale. As we increase the weight for the production cost function, the model continues to achieve the same, resulting in no major change in the network design.

\subsubsection{Impact of Economic Points}\label{sec:impacteconomic}
Next, we investigate the impact of economic points in the production cost function on the design of the network. Recall that the economic point $z_j^0$ is defined as $\beta K_j$, where $\beta$ is chosen as $0.5$. Increasing the value of $\beta$ shifts the economic point forward in the inverse S-shaped function. We vary $\beta$ from 0.5 to 0.75 and report the results in Table 9 (see Appendix E). From the result, we make the following observation.

\paragraph{Observation 1:} The total number of open facilities in the network either remains the same or reduces when $\beta$ is increased.\par
%
The reason for this is that the facilities operating under dis-economies of scale start operating under economies of scale with an increase in $\beta$, and their production also goes up, leading to the closure of a few facilities.

\section{Conclusions} \label{sec:conclusions}
In this paper, we have proposed an exact algorithm for solving a class of non-convex problems, where non-convexity arises due to the presence of one or more inverse S-shaped functions in the problem. In our approach, we have used an inner-approximation approach to iteratively refine the concave part of the inverse S-shaped function. Principles from bilevel optimization have been used to generate a lower bound formulation, which is a mixed integer convex program that can be handled using the cutting plane algorithm. The method guarantees global convergence over iterations. 
Using our approach, we solve 660 instances of facility location problems involving a linear fixed and transportation cost and inverse S-shaped production cost. The proposed algorithm is able to outperform the benchmark methods by converging to the global optimum for 95 percent of the instances within 3 hours of computation time. To our best knowledge, there does not exist general algorithms for solving non-convex problems with inverse S-shaped function globally, therefore we believe our algorithm will be very helpful in solving other application problems having similar structure.
\bibliography{mybib} 

\begin{thebibliography}{57}
\expandafter\ifx\csname natexlab\endcsname\relax\def\natexlab#1{#1}\fi
\providecommand{\url}[1]{\texttt{#1}}
\providecommand{\href}[2]{#2}
\providecommand{\path}[1]{#1}
\providecommand{\DOIprefix}{doi:}
\providecommand{\ArXivprefix}{arXiv:}
\providecommand{\URLprefix}{URL: }
\providecommand{\Pubmedprefix}{pmid:}
\providecommand{\doi}[1]{\href{http://dx.doi.org/#1}{\path{#1}}}
\providecommand{\Pubmed}[1]{\href{pmid:#1}{\path{#1}}}
\providecommand{\bibinfo}[2]{#2}
\ifx\xfnm\relax \def\xfnm[#1]{\unskip,\space#1}\fi
\bibitem[{A{\u{g}}ral{\i} \& Geunes(2009)}]{augrali2009solving}
\bibinfo{author}{A{\u{g}}ral{\i}, S.}, \& \bibinfo{author}{Geunes, J.}
  (\bibinfo{year}{2009}).
\newblock \bibinfo{title}{Solving knapsack problems with s-curve return
  functions}.
\newblock {\it \bibinfo{journal}{European Journal of Operational Research}\/},
  {\it \bibinfo{volume}{193}\/}, \bibinfo{pages}{605--615}.
\bibitem[{Alvarez \& Arias(2003)}]{alvarez2003diseconomies}
\bibinfo{author}{Alvarez, A.}, \& \bibinfo{author}{Arias, C.}
  (\bibinfo{year}{2003}).
\newblock \bibinfo{title}{Diseconomies of size with fixed managerial ability}.
\newblock {\it \bibinfo{journal}{American Journal of Agricultural
  Economics}\/},  {\it \bibinfo{volume}{85}\/}, \bibinfo{pages}{134--142}.
\bibitem[{Amiri(1997)}]{amiri1997solution}
\bibinfo{author}{Amiri, A.} (\bibinfo{year}{1997}).
\newblock \bibinfo{title}{Solution procedures for the service system design
  problem}.
\newblock {\it \bibinfo{journal}{Computers \& Operations Research}\/},  {\it
  \bibinfo{volume}{24}\/}, \bibinfo{pages}{49--60}.
\bibitem[{Belotti et~al.(2009)Belotti, Lee, Liberti, Margot \&
  W{\"a}chter}]{belotti2009branching}
\bibinfo{author}{Belotti, P.}, \bibinfo{author}{Lee, J.},
  \bibinfo{author}{Liberti, L.}, \bibinfo{author}{Margot, F.}, \&
  \bibinfo{author}{W{\"a}chter, A.} (\bibinfo{year}{2009}).
\newblock \bibinfo{title}{Branching and bounds tightening techniques for
  non-convex {MINLP}}.
\newblock {\it \bibinfo{journal}{Optimization Methods \& Software}\/},  {\it
  \bibinfo{volume}{24}\/}, \bibinfo{pages}{597--634}.
\bibitem[{Benjaafar et~al.(2004)Benjaafar, ElHafsi \&
  Vericourt}]{benjaafar2004demand}
\bibinfo{author}{Benjaafar, S.}, \bibinfo{author}{ElHafsi, M.}, \&
  \bibinfo{author}{Vericourt, F.~d.} (\bibinfo{year}{2004}).
\newblock \bibinfo{title}{Demand allocation in multiple-product,
  multiple-facility, make-to-stock systems}.
\newblock {\it \bibinfo{journal}{Management {S}cience}\/},  {\it
  \bibinfo{volume}{50}\/}, \bibinfo{pages}{1431--1448}.
\bibitem[{Bhatt et~al.(2021)Bhatt, Jayaswal, Sinha \&
  Vidyarthi}]{bhatt2021alternate}
\bibinfo{author}{Bhatt, S.~D.}, \bibinfo{author}{Jayaswal, S.},
  \bibinfo{author}{Sinha, A.}, \& \bibinfo{author}{Vidyarthi, N.}
  (\bibinfo{year}{2021}).
\newblock \bibinfo{title}{Alternate second order conic program reformulations
  for hub location under stochastic demand and congestion}.
\newblock {\it \bibinfo{journal}{Annals of Operations Research}\/},  (pp.
  \bibinfo{pages}{1--47}).
\bibitem[{Cohen \& Moon(1991)}]{cohen1991integrated}
\bibinfo{author}{Cohen, M.~A.}, \& \bibinfo{author}{Moon, S.}
  (\bibinfo{year}{1991}).
\newblock \bibinfo{title}{An integrated plant loading model with economies of
  scale and scope}.
\newblock {\it \bibinfo{journal}{European Journal of Operational Research}\/},
  {\it \bibinfo{volume}{50}\/}, \bibinfo{pages}{266--279}.
\bibitem[{Dasci \& Verter(2001)}]{dasci2001plant}
\bibinfo{author}{Dasci, A.}, \& \bibinfo{author}{Verter, V.}
  (\bibinfo{year}{2001}).
\newblock \bibinfo{title}{The plant location and technology acquisition
  problem}.
\newblock {\it \bibinfo{journal}{IIE {Transactions}}\/},  {\it
  \bibinfo{volume}{33}\/}, \bibinfo{pages}{963--974}.
\bibitem[{Daskin et~al.(2002)Daskin, Coullard \& Shen}]{daskin2002inventory}
\bibinfo{author}{Daskin, M.~S.}, \bibinfo{author}{Coullard, C.~R.}, \&
  \bibinfo{author}{Shen, Z.-J.~M.} (\bibinfo{year}{2002}).
\newblock \bibinfo{title}{An inventory-location model: Formulation, solution
  algorithm and computational results}.
\newblock {\it \bibinfo{journal}{Annals of {O}perations {R}esearch}\/},  {\it
  \bibinfo{volume}{110}\/}, \bibinfo{pages}{83--106}.
\bibitem[{Desrochers et~al.(1995)Desrochers, Marcotte \&
  Stan}]{desrochers1995congested}
\bibinfo{author}{Desrochers, M.}, \bibinfo{author}{Marcotte, P.}, \&
  \bibinfo{author}{Stan, M.} (\bibinfo{year}{1995}).
\newblock \bibinfo{title}{The congested facility location problem}.
\newblock {\it \bibinfo{journal}{Location Science}\/},  {\it
  \bibinfo{volume}{3}\/}, \bibinfo{pages}{9--23}.
\bibitem[{Dupont(2008)}]{dupont2008branch}
\bibinfo{author}{Dupont, L.} (\bibinfo{year}{2008}).
\newblock \bibinfo{title}{Branch and bound algorithm for a facility location
  problem with concave site dependent costs}.
\newblock {\it \bibinfo{journal}{International journal of {P}roduction
  {E}conomics}\/},  {\it \bibinfo{volume}{112}\/}, \bibinfo{pages}{245--254}.
\bibitem[{Elhedhli(2005)}]{elhedhli2005exact}
\bibinfo{author}{Elhedhli, S.} (\bibinfo{year}{2005}).
\newblock \bibinfo{title}{Exact solution of a class of nonlinear knapsack
  problems}.
\newblock {\it \bibinfo{journal}{Operations {R}esearch {L}etters}\/},  {\it
  \bibinfo{volume}{33}\/}, \bibinfo{pages}{615--624}.
\bibitem[{Elhedhli(2006)}]{elhedhli2006service}
\bibinfo{author}{Elhedhli, S.} (\bibinfo{year}{2006}).
\newblock \bibinfo{title}{Service system design with immobile servers,
  stochastic demand, and congestion}.
\newblock {\it \bibinfo{journal}{Manufacturing \& Service Operations
  Management}\/},  {\it \bibinfo{volume}{8}\/}, \bibinfo{pages}{92--97}.
\bibitem[{Elhedhli et~al.(2018)Elhedhli, Wang \& Saif}]{elhedhli2018service}
\bibinfo{author}{Elhedhli, S.}, \bibinfo{author}{Wang, Y.}, \&
  \bibinfo{author}{Saif, A.} (\bibinfo{year}{2018}).
\newblock \bibinfo{title}{Service system design with immobile servers,
  stochastic demand and concave-cost capacity selection}.
\newblock {\it \bibinfo{journal}{Computers \& Operations Research}\/},  {\it
  \bibinfo{volume}{94}\/}, \bibinfo{pages}{65--75}.
\bibitem[{Falk \& Soland(1969)}]{falk1969algorithm}
\bibinfo{author}{Falk, J.~E.}, \& \bibinfo{author}{Soland, R.~M.}
  (\bibinfo{year}{1969}).
\newblock \bibinfo{title}{An algorithm for separable nonconvex programming
  problems}.
\newblock {\it \bibinfo{journal}{Management {S}cience}\/},  {\it
  \bibinfo{volume}{15}\/}, \bibinfo{pages}{550--569}.
\bibitem[{Farahani et~al.(2015)Farahani, Rashidi~Bajgan, Fahimnia \&
  Kaviani}]{farahani2015location}
\bibinfo{author}{Farahani, R.~Z.}, \bibinfo{author}{Rashidi~Bajgan, H.},
  \bibinfo{author}{Fahimnia, B.}, \& \bibinfo{author}{Kaviani, M.}
  (\bibinfo{year}{2015}).
\newblock \bibinfo{title}{Location-inventory problem in supply chains: a
  modelling review}.
\newblock {\it \bibinfo{journal}{International Journal of Production
  Research}\/},  {\it \bibinfo{volume}{53}\/}, \bibinfo{pages}{3769--3788}.
\bibitem[{Feldman et~al.(1966)Feldman, Lehrer \& Ray}]{feldman1966warehouse}
\bibinfo{author}{Feldman, E.}, \bibinfo{author}{Lehrer, F.}, \&
  \bibinfo{author}{Ray, T.} (\bibinfo{year}{1966}).
\newblock \bibinfo{title}{Warehouse location under continuous economies of
  scale}.
\newblock {\it \bibinfo{journal}{Management Science}\/},  {\it
  \bibinfo{volume}{12}\/}, \bibinfo{pages}{670--684}.
\bibitem[{Fontes \& Gon{\c{c}}alves(2007)}]{fontes2007heuristic}
\bibinfo{author}{Fontes, D.~B.}, \& \bibinfo{author}{Gon{\c{c}}alves, J.~F.}
  (\bibinfo{year}{2007}).
\newblock \bibinfo{title}{Heuristic solutions for general concave minimum cost
  network flow problems}.
\newblock {\it \bibinfo{journal}{Networks: An International Journal}\/},  {\it
  \bibinfo{volume}{50}\/}, \bibinfo{pages}{67--76}.
\bibitem[{Guisewite \& Pardalos(1990)}]{guisewite1990minimum}
\bibinfo{author}{Guisewite, G.~M.}, \& \bibinfo{author}{Pardalos, P.~M.}
  (\bibinfo{year}{1990}).
\newblock \bibinfo{title}{Minimum concave-cost network flow problems:
  Applications, complexity, and algorithms}.
\newblock {\it \bibinfo{journal}{Annals of Operations Research}\/},  {\it
  \bibinfo{volume}{25}\/}, \bibinfo{pages}{75--99}.
\bibitem[{Hajiaghayi et~al.(2003)Hajiaghayi, Mahdian \&
  Mirrokni}]{hajiaghayi2003facility}
\bibinfo{author}{Hajiaghayi, M.~T.}, \bibinfo{author}{Mahdian, M.}, \&
  \bibinfo{author}{Mirrokni, V.~S.} (\bibinfo{year}{2003}).
\newblock \bibinfo{title}{The facility location problem with general cost
  functions}.
\newblock {\it \bibinfo{journal}{Networks: An International Journal}\/},  {\it
  \bibinfo{volume}{42}\/}, \bibinfo{pages}{42--47}.
\bibitem[{Holmberg(1994)}]{holmberg1994solving}
\bibinfo{author}{Holmberg, K.} (\bibinfo{year}{1994}).
\newblock \bibinfo{title}{Solving the staircase cost facility location problem
  with decomposition and piecewise linearization}.
\newblock {\it \bibinfo{journal}{European Journal of Operational Research}\/},
  {\it \bibinfo{volume}{75}\/}, \bibinfo{pages}{41--61}.
\bibitem[{Holmberg(1995)}]{holmberg1995efficient}
\bibinfo{author}{Holmberg, K.} (\bibinfo{year}{1995}).
\newblock \bibinfo{title}{Efficient decomposition and linearization methods for
  the stochastic transportation problem}.
\newblock {\it \bibinfo{journal}{Computational Optimization and
  Applications}\/},  {\it \bibinfo{volume}{4}\/}, \bibinfo{pages}{293--316}.
\bibitem[{Holmberg et~al.(1999)Holmberg, R{\"o}nnqvist \&
  Yuan}]{holmberg1999exact}
\bibinfo{author}{Holmberg, K.}, \bibinfo{author}{R{\"o}nnqvist, M.}, \&
  \bibinfo{author}{Yuan, D.} (\bibinfo{year}{1999}).
\newblock \bibinfo{title}{An exact algorithm for the capacitated facility
  location problems with single sourcing}.
\newblock {\it \bibinfo{journal}{European Journal of Operational Research}\/},
  {\it \bibinfo{volume}{113}\/}, \bibinfo{pages}{544--559}.
\bibitem[{Holmberg \& Tuy(1999)}]{holmberg1999production}
\bibinfo{author}{Holmberg, K.}, \& \bibinfo{author}{Tuy, H.}
  (\bibinfo{year}{1999}).
\newblock \bibinfo{title}{A production-transportation problem with stochastic
  demand and concave production costs}.
\newblock {\it \bibinfo{journal}{Mathematical {P}rogramming}\/},  {\it
  \bibinfo{volume}{85}\/}, \bibinfo{pages}{157--179}.
\bibitem[{Jayaswal et~al.(2017)Jayaswal, Vidyarthi \&
  Das}]{jayaswal2017cutting}
\bibinfo{author}{Jayaswal, S.}, \bibinfo{author}{Vidyarthi, N.}, \&
  \bibinfo{author}{Das, S.} (\bibinfo{year}{2017}).
\newblock \bibinfo{title}{A cutting plane approach to combinatorial bandwidth
  packing problem with queuing delays}.
\newblock {\it \bibinfo{journal}{Optimization Letters}\/},  {\it
  \bibinfo{volume}{11}\/}, \bibinfo{pages}{225--239}.
\bibitem[{Jeet et~al.(2009)Jeet, Kutanoglu \& Partani}]{jeet2009logistics}
\bibinfo{author}{Jeet, V.}, \bibinfo{author}{Kutanoglu, E.}, \&
  \bibinfo{author}{Partani, A.} (\bibinfo{year}{2009}).
\newblock \bibinfo{title}{Logistics network design with inventory stocking for
  low-demand parts: Modeling and optimization}.
\newblock {\it \bibinfo{journal}{IIE Transactions}\/},  {\it
  \bibinfo{volume}{41}\/}, \bibinfo{pages}{389--407}.
\bibitem[{Kelley(1960)}]{kelley1960cutting}
\bibinfo{author}{Kelley, J.~E., Jr} (\bibinfo{year}{1960}).
\newblock \bibinfo{title}{The cutting-plane method for solving convex
  programs}.
\newblock {\it \bibinfo{journal}{Journal of the Society for Industrial and
  Applied Mathematics}\/},  {\it \bibinfo{volume}{8}\/},
  \bibinfo{pages}{703--712}.
\bibitem[{Klincewicz et~al.(1990)Klincewicz, Luss \&
  Pilcher}]{klincewicz1990fleet}
\bibinfo{author}{Klincewicz, J.~G.}, \bibinfo{author}{Luss, H.}, \&
  \bibinfo{author}{Pilcher, M.~G.} (\bibinfo{year}{1990}).
\newblock \bibinfo{title}{Fleet size planning when outside carrier services are
  available}.
\newblock {\it \bibinfo{journal}{Transportation Science}\/},  {\it
  \bibinfo{volume}{24}\/}, \bibinfo{pages}{169--182}.
\bibitem[{Kuno \& Utsunomiya(2000)}]{kuno2000lagrangian}
\bibinfo{author}{Kuno, T.}, \& \bibinfo{author}{Utsunomiya, T.}
  (\bibinfo{year}{2000}).
\newblock \bibinfo{title}{A lagrangian based branch-and-bound algorithm for
  production-transportation problems}.
\newblock {\it \bibinfo{journal}{Journal of Global Optimization}\/},  {\it
  \bibinfo{volume}{18}\/}, \bibinfo{pages}{59--73}.
\bibitem[{Lee \& Grossmann(2001)}]{lee2001global}
\bibinfo{author}{Lee, S.}, \& \bibinfo{author}{Grossmann, I.~E.}
  (\bibinfo{year}{2001}).
\newblock \bibinfo{title}{A global optimization algorithm for nonconvex
  generalized disjunctive programming and applications to process systems}.
\newblock {\it \bibinfo{journal}{Computers \& Chemical Engineering}\/},  {\it
  \bibinfo{volume}{25}\/}, \bibinfo{pages}{1675--1697}.
\bibitem[{Li et~al.(2021)Li, Lin \& Shu}]{li2021location}
\bibinfo{author}{Li, Y.}, \bibinfo{author}{Lin, Y.}, \& \bibinfo{author}{Shu,
  J.} (\bibinfo{year}{2021}).
\newblock \bibinfo{title}{Location and two-echelon inventory network design
  with economies and diseconomies of scale in facility operating costs}.
\newblock {\it \bibinfo{journal}{Computers \& Operations Research}\/},  {\it
  \bibinfo{volume}{133}\/}, \bibinfo{pages}{105347}.
\bibitem[{Lin et~al.(2006)Lin, Nozick \& Turnquist}]{lin2006strategic}
\bibinfo{author}{Lin, J.-R.}, \bibinfo{author}{Nozick, L.~K.}, \&
  \bibinfo{author}{Turnquist, M.~A.} (\bibinfo{year}{2006}).
\newblock \bibinfo{title}{Strategic design of distribution systems with
  economies of scale in transportation}.
\newblock {\it \bibinfo{journal}{Annals of Operations Research}\/},  {\it
  \bibinfo{volume}{144}\/}, \bibinfo{pages}{161--180}.
\bibitem[{Lu et~al.(2014)Lu, Gzara \& Elhedhli}]{lu2014facility}
\bibinfo{author}{Lu, D.}, \bibinfo{author}{Gzara, F.}, \&
  \bibinfo{author}{Elhedhli, S.} (\bibinfo{year}{2014}).
\newblock \bibinfo{title}{Facility location with economies and diseconomies of
  scale: models and column generation heuristics}.
\newblock {\it \bibinfo{journal}{IIE {Transactions}}\/},  {\it
  \bibinfo{volume}{46}\/}, \bibinfo{pages}{585--600}.
\bibitem[{Marsten \& Morin(1978)}]{marsten1978hybrid}
\bibinfo{author}{Marsten, R.~E.}, \& \bibinfo{author}{Morin, T.~L.}
  (\bibinfo{year}{1978}).
\newblock \bibinfo{title}{A hybrid approach to discrete {M}athematical
  {P}rogramming}.
\newblock {\it \bibinfo{journal}{Mathematical {P}rogramming}\/},  {\it
  \bibinfo{volume}{14}\/}, \bibinfo{pages}{21--40}.
\bibitem[{McAfee \& McMillan(1995)}]{mcafee1995organizational}
\bibinfo{author}{McAfee, R.~P.}, \& \bibinfo{author}{McMillan, J.}
  (\bibinfo{year}{1995}).
\newblock \bibinfo{title}{Organizational diseconomies of scale}.
\newblock {\it \bibinfo{journal}{Journal of Economics \& Management
  Strategy}\/},  {\it \bibinfo{volume}{4}\/}, \bibinfo{pages}{399--426}.
\bibitem[{McCormick(1976)}]{mccormick1976computability}
\bibinfo{author}{McCormick, G.~P.} (\bibinfo{year}{1976}).
\newblock \bibinfo{title}{Computability of global solutions to factorable
  nonconvex programs: Part i—convex underestimating problems}.
\newblock {\it \bibinfo{journal}{Mathematical {P}rogramming}\/},  {\it
  \bibinfo{volume}{10}\/}, \bibinfo{pages}{147--175}.
\bibitem[{Mirchandani \& Jagannathan(1989)}]{mirchandani1989discrete}
\bibinfo{author}{Mirchandani, P.}, \& \bibinfo{author}{Jagannathan, R.}
  (\bibinfo{year}{1989}).
\newblock \bibinfo{title}{Discrete facility location with nonlinear
  diseconomies in fixed costs}.
\newblock {\it \bibinfo{journal}{Annals of {O}perations {R}esearch}\/},  {\it
  \bibinfo{volume}{18}\/}, \bibinfo{pages}{213--224}.
\bibitem[{Ni et~al.(2021)Ni, Shu, Song, Xu \& Zhang}]{ni2021branch}
\bibinfo{author}{Ni, W.}, \bibinfo{author}{Shu, J.}, \bibinfo{author}{Song,
  M.}, \bibinfo{author}{Xu, D.}, \& \bibinfo{author}{Zhang, K.}
  (\bibinfo{year}{2021}).
\newblock \bibinfo{title}{A branch-and-price algorithm for facility location
  with general facility cost functions}.
\newblock {\it \bibinfo{journal}{INFORMS Journal on Computing}\/},  {\it
  \bibinfo{volume}{33}\/}, \bibinfo{pages}{86--104}.
\bibitem[{Rockafellar(1970)}]{rockafellar1970convex}
\bibinfo{author}{Rockafellar, R.~T.} (\bibinfo{year}{1970}).
\newblock {\it \bibinfo{title}{Convex analysis}\/} volume~\bibinfo{volume}{36}.
\newblock \bibinfo{publisher}{Princeton university press}.
\bibitem[{Romeijn et~al.(2010)Romeijn, Sharkey, Shen \&
  Zhang}]{romeijn2010integrating}
\bibinfo{author}{Romeijn, H.~E.}, \bibinfo{author}{Sharkey, T.~C.},
  \bibinfo{author}{Shen, Z.-J.~M.}, \& \bibinfo{author}{Zhang, J.}
  (\bibinfo{year}{2010}).
\newblock \bibinfo{title}{Integrating facility location and production planning
  decisions}.
\newblock {\it \bibinfo{journal}{Networks: An International Journal}\/},  {\it
  \bibinfo{volume}{55}\/}, \bibinfo{pages}{78--89}.
\bibitem[{Ryoo \& Sahinidis(1995)}]{ryoo1995global}
\bibinfo{author}{Ryoo, H.~S.}, \& \bibinfo{author}{Sahinidis, N.~V.}
  (\bibinfo{year}{1995}).
\newblock \bibinfo{title}{Global optimization of nonconvex {NLP}s and {MINLP}s
  with applications in process design}.
\newblock {\it \bibinfo{journal}{Computers \& Chemical Engineering}\/},  {\it
  \bibinfo{volume}{19}\/}, \bibinfo{pages}{551--566}.
\bibitem[{Ryoo \& Sahinidis(1996)}]{ryoo1996branch}
\bibinfo{author}{Ryoo, H.~S.}, \& \bibinfo{author}{Sahinidis, N.~V.}
  (\bibinfo{year}{1996}).
\newblock \bibinfo{title}{A branch-and-reduce approach to global optimization}.
\newblock {\it \bibinfo{journal}{Journal of {G}lobal {O}ptimization}\/},  {\it
  \bibinfo{volume}{8}\/}, \bibinfo{pages}{107--138}.
\bibitem[{Saif(2016)}]{saif2016supply}
\bibinfo{author}{Saif, A.} (\bibinfo{year}{2016}).
\newblock \bibinfo{title}{Supply chain network design with concave costs:
  Theory and applications}.
\newblock {\it \bibinfo{journal}{University of Waterloo}\/},  (pp.
  \bibinfo{pages}{1--120}).
\bibitem[{Saif \& Elhedhli(2016)}]{saif2016lagrangian}
\bibinfo{author}{Saif, A.}, \& \bibinfo{author}{Elhedhli, S.}
  (\bibinfo{year}{2016}).
\newblock \bibinfo{title}{A lagrangian heuristic for concave cost facility
  location problems: the plant location and technology acquisition problem}.
\newblock {\it \bibinfo{journal}{Optimization Letters}\/},  {\it
  \bibinfo{volume}{10}\/}, \bibinfo{pages}{1087--1100}.
\bibitem[{Sharkey et~al.(2011)Sharkey, Geunes, Edwin~Romeijn \&
  Shen}]{sharkey2011exact}
\bibinfo{author}{Sharkey, T.~C.}, \bibinfo{author}{Geunes, J.},
  \bibinfo{author}{Edwin~Romeijn, H.}, \& \bibinfo{author}{Shen, Z.-J.~M.}
  (\bibinfo{year}{2011}).
\newblock \bibinfo{title}{Exact algorithms for integrated facility location and
  production planning problems}.
\newblock {\it \bibinfo{journal}{Naval Research Logistics (NRL)}\/},  {\it
  \bibinfo{volume}{58}\/}, \bibinfo{pages}{419--436}.
\bibitem[{Shectman \& Sahinidis(1998)}]{shectman1998finite}
\bibinfo{author}{Shectman, J.~P.}, \& \bibinfo{author}{Sahinidis, N.~V.}
  (\bibinfo{year}{1998}).
\newblock \bibinfo{title}{A finite algorithm for global minimization of
  separable concave programs}.
\newblock {\it \bibinfo{journal}{Journal of Global Optimization}\/},  {\it
  \bibinfo{volume}{12}\/}, \bibinfo{pages}{1--36}.
\bibitem[{Shen \& Qi(2007)}]{shen2007incorporating}
\bibinfo{author}{Shen, Z.-J.~M.}, \& \bibinfo{author}{Qi, L.}
  (\bibinfo{year}{2007}).
\newblock \bibinfo{title}{Incorporating inventory and routing costs in
  strategic location models}.
\newblock {\it \bibinfo{journal}{European Journal of Operational Research}\/},
  {\it \bibinfo{volume}{179}\/}, \bibinfo{pages}{372--389}.
\bibitem[{Shu et~al.(2005)Shu, Teo \& Shen}]{shu2005stochastic}
\bibinfo{author}{Shu, J.}, \bibinfo{author}{Teo, C.-P.}, \&
  \bibinfo{author}{Shen, Z.-J.~M.} (\bibinfo{year}{2005}).
\newblock \bibinfo{title}{Stochastic transportation-inventory network design
  problem}.
\newblock {\it \bibinfo{journal}{Operations Research}\/},  {\it
  \bibinfo{volume}{53}\/}, \bibinfo{pages}{48--60}.
\bibitem[{Sinha et~al.(2017)Sinha, Malo \& Deb}]{sinha2017review}
\bibinfo{author}{Sinha, A.}, \bibinfo{author}{Malo, P.}, \&
  \bibinfo{author}{Deb, K.} (\bibinfo{year}{2017}).
\newblock \bibinfo{title}{A review on bilevel optimization: from classical to
  evolutionary approaches and applications}.
\newblock {\it \bibinfo{journal}{IEEE Transactions on Evolutionary
  Computation}\/},  {\it \bibinfo{volume}{22}\/}, \bibinfo{pages}{276--295}.
\bibitem[{Soland(1971)}]{soland1971algorithm}
\bibinfo{author}{Soland, R.~M.} (\bibinfo{year}{1971}).
\newblock \bibinfo{title}{An algorithm for separable nonconvex programming
  problems ii: Nonconvex constraints}.
\newblock {\it \bibinfo{journal}{Management {S}cience}\/},  {\it
  \bibinfo{volume}{17}\/}, \bibinfo{pages}{759--773}.
\bibitem[{Soland(1974)}]{soland1974optimal}
\bibinfo{author}{Soland, R.~M.} (\bibinfo{year}{1974}).
\newblock \bibinfo{title}{Optimal facility location with concave costs}.
\newblock {\it \bibinfo{journal}{Operations Research}\/},  {\it
  \bibinfo{volume}{22}\/}, \bibinfo{pages}{373--382}.
\bibitem[{Solomon(1987)}]{solomon1987algorithms}
\bibinfo{author}{Solomon, M.~M.} (\bibinfo{year}{1987}).
\newblock \bibinfo{title}{Algorithms for the vehicle routing and scheduling
  problems with time window constraints}.
\newblock {\it \bibinfo{journal}{Operations {R}esearch}\/},  {\it
  \bibinfo{volume}{35}\/}, \bibinfo{pages}{254--265}.
\bibitem[{Taha(1973)}]{taha1973concave}
\bibinfo{author}{Taha, H.~A.} (\bibinfo{year}{1973}).
\newblock \bibinfo{title}{Concave minimization over a convex polyhedron}.
\newblock {\it \bibinfo{journal}{Naval Research Logistics Quarterly}\/},  {\it
  \bibinfo{volume}{20}\/}, \bibinfo{pages}{533--548}.
\bibitem[{Tawarmalani \& Sahinidis(2005)}]{tawarmalani2005polyhedral}
\bibinfo{author}{Tawarmalani, M.}, \& \bibinfo{author}{Sahinidis, N.~V.}
  (\bibinfo{year}{2005}).
\newblock \bibinfo{title}{A polyhedral branch-and-cut approach to global
  optimization}.
\newblock {\it \bibinfo{journal}{Mathematical {P}rogramming}\/},  {\it
  \bibinfo{volume}{103}\/}, \bibinfo{pages}{225--249}.
\bibitem[{Tuy(1964)}]{tuy1964concave}
\bibinfo{author}{Tuy, H.} (\bibinfo{year}{1964}).
\newblock \bibinfo{title}{Concave programming under linear constraints}.
\newblock {\it \bibinfo{journal}{Soviet Math.}\/},  {\it
  \bibinfo{volume}{5}\/}, \bibinfo{pages}{1437--1440}.
\bibitem[{Vidyarthi \& Jayaswal(2014)}]{vidyarthi2014efficient}
\bibinfo{author}{Vidyarthi, N.}, \& \bibinfo{author}{Jayaswal, S.}
  (\bibinfo{year}{2014}).
\newblock \bibinfo{title}{Efficient solution of a class of location--allocation
  problems with stochastic demand and congestion}.
\newblock {\it \bibinfo{journal}{Computers \& Operations Research}\/},  {\it
  \bibinfo{volume}{48}\/}, \bibinfo{pages}{20--30}.
\bibitem[{Zetina et~al.(2021)Zetina, Contreras \& Jayaswal}]{zetina2021exact}
\bibinfo{author}{Zetina, C.~A.}, \bibinfo{author}{Contreras, I.}, \&
  \bibinfo{author}{Jayaswal, S.} (\bibinfo{year}{2021}).
\newblock \bibinfo{title}{An exact algorithm for large-scale non-convex
  quadratic facility location}.
\newblock {\it \bibinfo{journal}{arXiv preprint arXiv:2107.09746}\/}, .

\end{thebibliography}
\newpage
\appendix
\hspace{3cm}\Huge \textbf{Online Supplements: Appendices}\\
\setcounter{page}{1}
\normalsize
\section{Theorems and Proofs}
\begin{theorem}\label{th:bigM}
$M_j^0=\phi_j^{\cap}(K_j)$ and $M_j^1=\phi_j^{\cup}(K_j)$ are acceptable values for large numbers $M_j^0$ and $M_j^1$.
\end{theorem}
\begin{proof}
From \eqref{eq:origConst9} and \eqref{eq:origConst10} , if $l_{j0}= 1$, then $w_j \le M_j^0, w_j \geq \phi_j^{\cap}(z_j) \implies w_j= \phi_j^{\cap}(z_j), w_j \le M_j^0$. if $l_{j0}= 0$, then $w_j \le 0, w_j \geq \phi_j^{\cap}(z_j)-M_j^0$, which indicates that we should choose $M_j^0,$ so that $w_j \le M_j^0 \; \forall \; z_j$. Since $z_j$ is bounded by $e_j \le z_j \le K_j$, the maximum possible value of $w_j$ is $\max\{\phi_j^{\cap}(z_j): e_j \leq z_j \leq K_j\}=\phi_j^{\cap}(K_j)$. Hence, $M_j^0=\phi_j^{\cap}(K_j)$ is acceptable.\par
Similarly, from \eqref{eq:origConst11} and \eqref{eq:origConst12} , if $l_{j1}= 1$, then $p_j \le M_j^1, p_j \geq \phi_j^{\cup}(z_j) \implies p_j= \phi_j^{\cup}(z_j), p_j \le M_j^1$. If $l_{j1}= 0$, then $p_j \le 0, p_j \geq \phi_j^{\cup}(z_j)-M_j^1$, which indicates that we should select $M_j^1,$ such that $p_j \le M_j^1 \; \forall z_j$. The maximum possible value of $p_j$ is $\max \{\phi_j^{\cup}(z_j): e_j \le z_j \le K_j\}=\phi_j^{\cup}(K_j)$. Hence, $M_j^1=\phi_j^{\cup}(K_j)$ is acceptable.
\end{proof}
\begin{theorem}\label{th:equal}
The IOA algorithm will terminate with optimal solution for \eqref{eq:origObj}-\eqref{eq:origConst5} if two consecutive iterations $c$ and $c+1$ in Algorithm~\ref{IOAalgo} lead to the same solution.
\end{theorem}
\begin{proof}
Let $x^{\tau+c},y^{\tau+c},q^{\tau+c},w^{\tau+c},p^{\tau+c}$ and $x^{\tau+c+1},y^{\tau+c+1},q^{\tau+c+1},w^{\tau+c+1},p^{\tau+c+1}$ be the solutions for the variables $x,y,z,w,p$ at iterations $c$ and $c+1$, respectively, for formulation \eqref{eq:origObj5}-\eqref{eq:origConst22}. It is important to note that $S^{c+1} = S^c \cup q^{\tau+c}$, which implies $\hat{\phi}_j^\cap(q^{\tau+c}|S^{c+1}) = \phi_j^\cap(q^{\tau+c}), j=1,\ldots, J$ at iteration $c+1$. At iteration $c+1$, $f(x^{\tau+c+1}, y^{\tau+c+1}) + \sum_{j=1}^{J} w_j^{\tau+c+1}+ \sum_{j=1}^{J} p_j^{\tau+c+1} \le f(x^{\tau+c+1}, y^{\tau+c+1}) + \sum_{j=1}^{J} \phi_j(z^{\tau+c+1})$. Since $x^{\tau+c+1}=x^{\tau+c},y^{\tau+c+1}=y^{\tau+c},q^{\tau+c+1}=q^{\tau+c},w^{\tau+c+1}=w^{\tau+c}, p^{\tau+c+1}=p^{\tau+c}$, and $\hat{\phi}_j^{\cap}(q^{\tau+c}|S^{c+1}) = \phi_j^{\cap}(q^{\tau+c})$, $f(x^{\tau+c+1}, y^{\tau+c+1}) + \sum_{j=1}^{J}w_j^{\tau+c+1}+\sum_{j=1}^{J}p_j^{\tau+c+1} = f(x^{\tau+c+1}, y^{\tau+c+1}) + \sum_{j=1}^{J} \phi_j(z^{\tau+c+1})$ (i.e. holds with an equality), which implies $x^{\tau+c},y^{\tau+c},q^{\tau+c},w^{\tau+c},p^{\tau+c}$ is the optimal solution for the variables $x,y,z,w,p$, respectively.
\end{proof}
\begin{theorem}
As the IOA algorithm moves from iteration $c$ to $c+1$, the lower bound for \eqref{eq:origObj}-\eqref{eq:origConst5} improves if the solution for Mod-$S^c$ (\eqref{eq:origObj5}-\eqref{eq:origConst22}) at iteration $c$ is not the optimum for \eqref{eq:origObj}-\eqref{eq:origConst5}.
\end{theorem}
\begin{proof}
Say $x^{\tau+c},y^{\tau+c},q^{\tau+c},w^{\tau+c},p^{\tau+c}$ be the solution for the variables $x,y,z,w,p$, respectively, at iteration $c$. Assume that this solution for formulation Mod-$S^c$ at iteration $c$ is not optimal for the original problem \eqref{eq:origObj}-\eqref{eq:origConst5}. Let $x^{\tau+c+1},y^{\tau+c+1},q^{\tau+c+1},w^{\tau+c+1},p^{\tau+c+1}$ be the solution for the variables $x,y,z,w,p$, respectively, at iteration $c+1$. From Theorem~\ref{th:equal} we can claim that the solution at iteration $c$ is not the same as the solution at iteration $c+1$.\par
If we add $\mu_{j(\tau+1)}=0, j= 1,\ldots, J$ in the linear program corresponding to $\hat{\phi_j^{\cap}}(z_j|S^{c+1})$, it becomes equivalent to the linear program corresponding to $\hat{\phi}_j^{\cap}(z_j|S^{c})$. Hence, the linear program corresponding to $\hat{\phi_j^{\cap}}(z_j|S^{c})$ is a relaxation of linear program corresponding to $\hat{\phi_j^{\cap}}(z_j|S^{c+1})$. So, for any given $z_j$, $\hat{\phi_j^{\cap}}(z_j|S^{c}) \le \hat{\phi_j^{\cap}}(z_j|S^{c+1})$. This implies that for \eqref{eq:origObj3}-\eqref{eq:origConst12}, $f(x^{\tau+c},y^{\tau+c})+\sum_{j=1}^{J}w_j^{\tau+c} + \sum_{j=1}^{J}p_j^{\tau+c} \le f(x^{\tau+c+1},y^{\tau+c+1})+\sum_{j=1}^{J}w_j^{\tau+c+1} + \sum_{j=1}^{J} p_j^{\tau+c+1}$. The equality will not hold due to strict concavity of $\phi_j^{\cap}(z_j)$ and the solutions being different ($x^{\tau+c} \neq x^{\tau+c+1}, y^{\tau+c} \neq y^{\tau+c+1}, w^{\tau+c} \neq w^{\tau+c+1}, p^{\tau+c} \neq p^{\tau+c+1}$). Therefore, $f(x^{\tau+c},y^{\tau+c})+\sum_{j=1}^{J}w_j^{\tau+c} + \sum_{j=1}^{J}p_j^{\tau+c} < f(x^{\tau+c+1},y^{\tau+c+1})+\sum_{j=1}^{J}w_j^{\tau+c+1} + \sum_{j=1}^{J} p_j^{\tau+c+1}$, which proves that the lower bound for \eqref{eq:origObj}-\eqref{eq:origConst5} strictly improves in the succeeding iteration. 
\end{proof}

\section{An Illustrative Example}\label{sec:NumEx}
To illustrate the algorithm, we consider a small-size numerical example:
\begin{align}
	\min_{x} & \; \;  x_1^{3} - 15x_1^2 + 75x_1 - 30x_2 \label{eq:exorigObj}\\
	\st & \;\;\; -9x_1 + 5x_2 \leq 9 \label{eq:exorigcons1} \\
	& \;\;\; x_1 - 6x_2 \leq 6 \label{eq:exorigcons2}\\
	& \;\;\; 3x_1 + x_2 \leq 9 \label{eq:exorigcons3}\\
	& x \in X = \{1 \leq x_j \leq 7, j = 1,2 \} \label{eq:exorigcons4}
\end{align}
where $\phi(x_1) =x_1^{3} - 15x_1^2 + 75x_1$ is an inverted S-shaped function.\\
Iteration 1: We split the inverse S-shaped function $x_1^{3} - 15 x_1^2 + 75 x_1$ into separate concave ($\phi^\cap(x_1)$) and convex ($\phi^\cup(x_1)$) functions as follows:
	\begin{equation}
	\phi^{\cap}(x_1) = 
	\begin{cases}
	x_1^{3} - 15x_1^2 + 75x_1 & \forall \; x_1\in [1,5]\\
	125 \ & \forall \; x_1 \in [5,7] \label{eq:exfundef1}
	\end{cases}       
	\end{equation}
	\begin{equation}
	\phi^{\cup}(x_1) =
	\begin{cases}
	125 & \forall \; x_1 \in [1,5]\\
	x_1^{3} - 15x_1^2 + 75x_1 & \forall \; x_1\in [5,7] \label{eq:exfundef2}
	\end{cases}       
	\end{equation}
Hence, the modified problem would be
\begin{align*}
	\min_{x} & \; \;  l_1 \phi^{\cap}(x_1)+l_2 \phi^{\cup}(x_1) - 30 x_2 \\
	\st &  -9x_1 + 5x_2 \leq 9 \\
	& \;\;\; x_1 - 6x_2 \leq 6 \\
	& \;\;\; 3x_1 + x_2 \leq 9 \\
	& \;\;\;  x \in X = \{1 \leq x_j \leq 7, j = 1,2 \}\\
	& \;\;\; l_1+l_2 \leq 1  \\
	& \;\;\; 5l_2\leq x_1 \leq 7l_2+5 l_{1}   \\
	& \;\;\; l_1,l_2 \in \{0,1\} 
\end{align*}
Next, we separate binary variables $l_1$ and $l_2$ from the non-convex terms $l_1 \phi^{\cap}(x_1)$ and $l_2 \phi^{\cup}(x_1)$, respectively using bigM method.
\begin{align}
	\min_{x} & \; \;  w + p - 30 x_2 \label{eq:exObj}\\
	\st &  -9x_1 + 5x_2 \leq 9 \label{eq:exConst1}\\
	& \;\;\; x_1 - 6x_2 \leq 6 \label{eq:exConst2} \\
	& \;\;\; 3x_1 + x_2 \leq 9 \label{eq:exConst3}\\
	& \;\;\;  x \in X = \{1 \leq x_j \leq 7,\; j = 1,2 \} \label{eq:exConst4}\\
	& \;\;\; l_1+l_2 \leq 1  \label{eq:exConst5}\\
	& \;\;\; 5l_2\leq x_1 \leq 7l_2+5 l_1  \label{eq:exConst6} \\
	& \;\;\; M^0l_1 \geq w  \label{eq:exConst7}\\
	& \;\;\; M^0l_1+ \phi^{\cap}(x_1) - M^0 \leq w  \label{eq:exConst8}\\
	& \;\;\; M^1l_2 \geq p, \; j=1,\dots,J \label{eq:exConst9}\\
	& \;\;\; M^1 l_2 + \phi^{\cup}(x_1) - M^1 \leq p \label{eq:exConst10}\\
	& \;\;\; l_1,l_2 \in \{0,1\}  \label{eq:exConst11}
\end{align}
We select $M^0= 125$ and $M^1=133$ as per theorem proposed in Theorem~\ref{th:bigM}.
Thereafter, we relax the concave constraints with the inner-approximation generated using two starting points, $x_1 \in \{1,7\}$, which gives us the following bilevel program.
Let $g(x_1) = \phi^{\cap}(x_1) $, then $g(1) = 1$, $g(7) = 125$.
\begin{align}
\min_{x} & \; \; \phi(x) = w + p - 30 x_2 \label{eq:exobj2}\\
	\st  \;\;\;
	& \mu \in \argmax_{\mu} \left\{ \mu_1 + 125 \mu_2 : \mu_1 + \mu_2 = 1, \mu_1 + 7\mu_2 = x_1,-\mu_1 \leq 0 -\mu_2 \leq 0 \right\} \label{eq:exConst12}\\
	&  t_1 \ge \mu_1 + 125 \mu_2 \label{eq:exConst13}\\
	& M^0l_1+ t_1 - M^0 \leq w \label{eq:exConst14}\\
	& \eqref{eq:exConst1}- \eqref{eq:exConst7},\eqref{eq:exConst9}-\eqref{eq:exConst11} \notag
\end{align}
Let $\gamma_1,\gamma_2,\gamma_3, \textrm{and}, \gamma_4$ be the Lagrange multipliers for the constraints in \eqref{eq:exConst12}, then the lower level problem in \eqref{eq:exConst12} can be converted to single level problem using following KKT conditions:
\begin{align}
& -1 + \gamma_1 + \gamma_2 - \gamma_3 = 0 \label{eq:excons15}\\
& -125 + \gamma_1 + 7 \gamma_2 - \gamma_4 = 0 \label{eq:excons16}\\
&-\mu_1 \gamma_3 = 0 \label{eq:excons17comp1}\\
&-\mu_2 \gamma_4 = 0 \label{eq:excons17comp2}\\
&\mu_1, \mu_2 , \gamma_3, \gamma_4 \geq 0 \label{eq:excons18}\\ 
&\gamma_1, \gamma_2 - \textrm{unrestricted} \label{eq:excons19}
\end{align}
We linearize complimentary slackness conditions \eqref{eq:excons17comp1}- \eqref{eq:excons17comp2} using the BigM values.
\begin{align}
&\mu_1 \leq M Z_1 \label{eq:excons20}\\
&\gamma_3 \leq M (1 - Z_1) \label{eq:excons21}\\
&\mu_2 \leq M Z_2 \label{eq:excons22}\\
&\gamma_4 \leq M (1 - Z_2) \label{eq:excons23}\\
&Z_1 , Z_2 \in \{0,1\} \label{eq:excons24}
\end{align}
The relaxed model for the original problem (\eqref{eq:exorigObj}-\eqref{eq:exorigcons4}) is given below as a mixed integer convex program. 
\begin{align*}
[EX1-1] \min & w + p - 30 x_2 \\
\st &\;\;\; \eqref{eq:exConst1}- \eqref{eq:exConst7},\eqref{eq:exConst9}-\eqref{eq:exConst11},\eqref{eq:exConst13}-\eqref{eq:exConst14}, \eqref{eq:excons15}-\eqref{eq:excons16},\eqref{eq:excons18}-\eqref{eq:excons24}
\end{align*}
The above formulation can be solved globally using a cutting plane algorithm to reach to the following solution, $x_1 = 1.5, x_2 = 4.5, \textrm{ objective value } = -108.214$.
Hence, the lower bound is -108.214 and the upper bound is -52.875. The optimality gap ($\%$) is 51.14 $\%$. Additional iterations would lead to improvement of the bounds and converge to global optimum after a few iterations. The global optimum for this problem is -52.875.

\section{Implementation of the Cutting Plane Method}\label{sec:cuttingPlane}
The replacement of concave function with it's inner approximation and further linearization of complementary slackness conditions convert the non-convex problem (\eqref{eq:origObj3}-\eqref{eq:origConst12}) to CMIP (\eqref{eq:origObj5}-\eqref{eq:origConst18} , \eqref{eq:origConst20} -\eqref{eq:lin3}), which can be solved globally using \cite{kelley1960cutting}'s cutting plane technique. 
In order to linearize the convex non-linear function $\phi_{j}^{\cup}(z_j)$, we define additional non-negative variable $s_j, \; \forall j$.
$\phi_j^{\cup}(z_j)$ can be outer approximated with the addition of following tangents at points $(z_j^h)_{h \in H_t}:\phi_{j}^{\cup}(z_j) \approx \phi_{j}^{\cup}(z_j^h) + \phi_{j}^{\cup'}(z_j^h)(z_j-z_j^h)$. Hence, constraint \eqref{eq:origConst12} can be approximated as follows:
\begin{align}
	& M_j^1 l_{j1} + s_j - M^1_j \leq p_j   \;\;\forall\; j=1,\dots,J \label{eq:origConst23} \\
	& s_j \geq \phi_j^{\cup}(z_j^h) + \phi_j^{\cup'}(z_j^h)(z_j-z_j^h), \;\;  \; \; h \in H_t; \;\; j=1,\dots,J \label{eq:origConst24}
\end{align}
The above substitution will lead to the following linear mixed integer program (MIP).
\begin{align}
	\textbf{[Mod-$S^c$-Linear]} \quad \min_{x} & \;  f(x,y) + \sum_{j=1}^{J}w_j + \sum_{j=1}^{J}p_j \label{eq:origObjMILP}\ \\
	\st & \; \; \; \eqref{eq:origConst} - \;\eqref{eq:origConst5},\eqref{eq:origConst6} - \;\eqref{eq:origConst8},\eqref{eq:origConst9}, \eqref{eq:origConst11},\eqref{eq:origConst13} \label{eq:origConst25}\\
	& \eqref{eq:origConst15} - \;\eqref{eq:origConst18},\eqref{eq:origConst20}-\;\eqref{eq:lin3},\eqref{eq:origConst23} - \eqref{eq:origConst24} \label{eq:origConst26} 
\end{align}
For any subset of points $\{s_j^h\}_{ h \in H_t \subset H}$, [Mod-$S^c$-Linear] is a relaxed version of the original problem (\eqref{eq:origObj5}-\eqref{eq:origConst18}, \eqref{eq:origConst20} -\eqref{eq:lin3}). Hence, $\sum_{j=1}^{J} p_j(H_t) $ will be the lower bound for the cutting plane algorithm, since $p_j(H_t) \leq p_j(H)$. Similarly, for any subset of points $\{s_j^h\}_{ h \in H_t \subset H},\;$ the solution $(z_j^h)$ is also a feasible solution. Hence $\sum_{j=1}^{J} \phi_j^{\cup}(z_j^h)$ will be the upper bound for the cutting plane algorithm.
To solve Mod-$S^c$ more efficiently, the cutting plane algorithm starts with \textit{a priori} set of cuts on carefully selected points \citep{elhedhli2005exact,jayaswal2017cutting}. A psuedo-code of cutting plane algorithm is provided below.\par
\begin{algorithm}
	\caption*{Cutting Plane Algorithm \citep{kelley1960cutting} }\label{CuttingPlaneAlgo}
	\begin{algorithmic}[1]
		\State \textit{Begin}
		\State $UB_\mathcal{C} \gets +\infty, LB_\mathcal{C} \gets -\infty, t \gets \textit{1}$
		\State $ \text{Choose an initial set of $\beta$ points } H_t = \{ v^{1}, v^{2}, \ldots, v^{\beta} \} $ as per apriori sets \citep{elhedhli2005exact}. 
		\While {(($UB_\mathcal{C}-LB_\mathcal{C})/LB_\mathcal{C} > \epsilon \; ) $ begin} 
		\State Solve Mod-$S^c$-Linear (\eqref{eq:origObjMILP}-\eqref{eq:origConst26}) with an MIP solver.
		\State Let the optimal solution for Mod-$S^c$-Linear be  $v^{\beta+t},p^{\beta +t},l^{\beta +t},\zeta^{\beta +t} \text{ for the variables $z,p,l,\zeta$ respectively}$
		\State $LB_C \gets f(x^{\beta +t},y^{\beta +t}) + \sum_{j=1}^{J}l_{j0}^{\beta +t} \zeta_j^{\beta +t} +\sum_{j=1}^{J}p^{\beta +t}_j $
		\State $ UB_C \gets f(x^{\beta +t},y^{\beta +t}) + \sum_{j=1}^{J}l_{j0}^{\beta +t} \zeta_j^{\beta +t}+ \sum_{j=1}^{J} l_{j1}^{\beta+t}\phi_j^{\cup}(v_j^{\beta+t})$
		\State $  \mathcal{H}_{t+1} \gets  \mathcal{H}_{t} \; \cup \;v^{\beta +t}  $
		\State $t \gets t +1$
		\EndWhile
		\State \textit{End}
	\end{algorithmic}
\end{algorithm}
\section{Data Synthesis Strategy}
Here, we represent the data synthesis strategy for facility $j \in N$ and customer $i \in M$. Note that we have broken Set 1 into two parts because of significant differences in the problems studied in this set.
\begin{itemize}
    \item [Set 1a] Problem 1-12 (p1-p12): $n=10,m=50,L \sim Uniform(10,200), d_i \sim Uniform(10,50), F_j \sim Uniform(300,700), K_j \sim Uniform(100,500), c_{ij}=\rho e_{ij},\rho=4, i=1,\ldots,m,j=1,\ldots,n.$
    \item [Set 1b] Problem 13-24 (p13-p24): $n=20,m=50, L \sim Uniform(10,200), d_i \sim Uniform(30,80),F_j \sim Uniform(300,700), K_j \sim Uniform(100,500),c_{ij}=\rho e_{ij},\rho=4, i=1,\ldots,m,j=1,\ldots,n.$ 
    \item [Set 2 ] Problem 25-40 (p25-p40): $n=30,m=150, L \sim Uniform(10,300), d_i \sim Uniform(10,50),F_j \sim [300,700], K_j \sim Uniform(200,600),c_{ij}=\rho_1(dw_j +e_{ij}+dw_{i}) +\rho_2 a_j, i=1,\ldots,m,j=1,\ldots,n.$ The transportation costs ($c_{ij}, j \in N$) are based on a vehicle routing problem modified from \citep{klincewicz1990fleet}. Here $dw_{j}$ is the distance between depot and the location of customer $j$ assuming depot is at $(0,0)$. $e_{ij}$ is the distance between location of $i$ and $j$. $\rho_1=4,\rho_2=2$.
    \item [Set 3 ] Problem 41-55 (p41-p55): $n \in \{10,20,30\},m \in \{70,80,90,100\}$. The coordinates of the locations and the demand is generated from random, random-clustered, and clustered-clustered problems of \cite{solomon1987algorithms}. $F_j \sim Uniform(300,600),K_j \sim Uniform(100,500)$. For problems p41-p49, the transportation cost $c_{ij}$ is similar to p1-p24, whereas for p50-p55 $c_{ij}$ is similar to p25-p40.
\end{itemize}
\begin{table}[htbp]
  \centering
  \caption{Test instances}
    \begin{tabular}{cccccccccrrr}
    \toprule
    Instance & $n,m$ & $K/D$   & Instance & $n,m$ & K/D   & Instances & $n,m$ & K/D   & \multicolumn{1}{c}{Instance} & \multicolumn{1}{c}{$n,m$} & \multicolumn{1}{c}{K/D} \\
    \midrule
    p1-p4 & 10,50 & 1.74  & p29-p32 & 30,150 & 3.03  & p45   & 20,80 & 4.14  & \multicolumn{1}{c}{p52} & \multicolumn{1}{c}{10,100} & \multicolumn{1}{c}{1.6} \\
    p5-p8 & 10,50 & 1.37  & p33-p36 & 30,150 & 4.04  & p46   & 30,70 & 7.1   & \multicolumn{1}{c}{p53} & \multicolumn{1}{c}{20,100} & \multicolumn{1}{c}{3.37} \\
    p9-p12 & 10,50 & 2.06  & p37-p40 & 30,150 & 6.06  & p47   & 10,90 & 1.76  & \multicolumn{1}{c}{p54} & \multicolumn{1}{c}{10,100} & \multicolumn{1}{c}{1.52} \\
    p13-p16 & 20,50 & 2.77  & p41   & 10,90 & 2.12  & p48   & 20,80 & 4.06  & \multicolumn{1}{c}{p55} & \multicolumn{1}{c}{20,100} & \multicolumn{1}{c}{3.21} \\
    p17-p20 & 20,50 & 2.8   & p42   & 20,80 & 4.99  & p49   & 30,70 & 7.08  &       &       &  \\
    p21-p24 & 20,50 & 3.5   & p43   & 30,70 & 8.28  & p50   & 10,100 & 1.89  &       &       &  \\
    p25-p28 & 30,150 & 4.12  & p44   & 10,90 & 1.76  & p51   & 20,100 & 3.98  &       &       &  \\
    \bottomrule
    \end{tabular}%
  \label{tab:datacharac}%
\end{table}%

\section{Analysis of Impact of Different Cost Components and Economic Point}
\begin{table}[htbp]
  \centering
  \caption{Impact of the transportation cost for function type 2 on the design of the nework}
     \resizebox{.94\textwidth}{!}{
    \begin{tabular}{ccccccccccccccccccc}
    \toprule
       &    & \multicolumn{5}{c}{ $\alpha=1$ } &    & \multicolumn{5}{c}{$\alpha=10$} &    & \multicolumn{5}{c}{$\alpha=20$} \\
\cmidrule{3-7}\cmidrule{9-13}\cmidrule{15-19}    Data  &    & Obj/Gap (\%) & Time & $n_e$ & $n_d$ & $n_T$ &    & Obj/Gap (\%) & Time & $n_e$ & $n_d$ & $n_T$ &    & Obj/Gap (\%) & Time & $n_e$ & $n_d$ & $n_T$ \\
\midrule
    p1 &    & 14468.3 & 0.82 & 0  & 8  & 8  &    & 60793.4 & 0.60 & 4  & 6  & 10 &    & 110683 & 0.46 & 4  & 6  & 10 \\
    p2 &    & 13494.2 & 9.48 & 0  & 9  & 9  &    & 59479.4 & 0.50 & 4  & 6  & 10 &    & 109369 & 0.44 & 4  & 6  & 10 \\
    p3 &    & 15159.3 & 1.10 & 0  & 8  & 8  &    & 61479.4 & 0.77 & 4  & 6  & 10 &    & 111369 & 0.37 & 4  & 6  & 10 \\
    p4 &    & 16599.2 & 0.29 & 0  & 7  & 7  &    & 63479.2 & 0.30 & 4  & 6  & 10 &    & 113369 & 0.31 & 4  & 6  & 10 \\
    p5 &    & 15338 & 0.39 & 0  & 10 & 10 &    & 60647.5 & 0.18 & 1  & 9  & 10 &    & 109721 & 0.40 & 1  & 9  & 10 \\
    p6 &    & 14024.1 & 0.42 & 0  & 10 & 10 &    & 59333.5 & 0.18 & 1  & 9  & 10 &    & 108405 & 0.45 & 1  & 9  & 10 \\
    p7 &    & 16024.1 & 0.43 & 0  & 10 & 10 &    & 61333.5 & 0.20 & 1  & 9  & 10 &    & 110405 & 0.39 & 1  & 9  & 10 \\
    p8 &    & 18024 & 0.44 & 0  & 10 & 10 &    & 63333.7 & 0.23 & 1  & 9  & 10 &    & 112407 & 0.42 & 1  & 9  & 10 \\
    p9 &    & 13232.2 & 5.04 & 1  & 7  & 8  &    & 56676.8 & 0.75 & 5  & 5  & 10 &    & 104018 & 0.37 & 5  & 5  & 10 \\
    p10 &    & 12282 & 14.37 & 1  & 7  & 8  &    & 55362.8 & 0.64 & 5  & 5  & 10 &    & 102704 & 0.44 & 5  & 5  & 10 \\
    p11 &    & 13852.1 & 9.98 & 0  & 7  & 7  &    & 57362.7 & 0.71 & 5  & 5  & 10 &    & 104704 & 0.41 & 5  & 5  & 10 \\
    p12 &    & 15252.1 & 0.92 & 0  & 7  & 7  &    & 59362.7 & 0.72 & 5  & 5  & 10 &    & 106704 & 0.40 & 5  & 5  & 10 \\
    p13 &    & 16876.5 & 72.56 & 0  & 10 & 10 &    & 49996.3 & 107.89 & 13 & 5  & 18 &    & 81972.2 & 2.25 & 15 & 5  & 20 \\
    p14 &    & 15287.7 & 370.08 & 0  & 10 & 10 &    & 47310.4 & 26.42 & 14 & 5  & 19 &    & 79120.2 & 2.47 & 15 & 5  & 20 \\
    p15 &    & 17225.3 & 11.25 & 0  & 9  & 9  &    & 50850.9 & 34.55 & 11 & 6  & 17 &    & 83120.2 & 2.83 & 15 & 5  & 20 \\
    p16 &    & 19025.2 & 2.63 & 0  & 9  & 9  &    & 54250.8 & 31.65 & 11 & 6  & 17 &    & 87120.2 & 1.65 & 15 & 5  & 20 \\
    p17 &    & 16710 & 166.00 & 0  & 10 & 10 &    & 49993.5 & 107.17 & 13 & 5  & 18 &    & 82045.2 & 2.69 & 15 & 5  & 20 \\
    p18 &    & 15377.4 & 1037.00 & 0  & 11 & 11 &    & 47368 & 33.07 & 13 & 5  & 18 &    & 79193.2 & 2.66 & 15 & 5  & 20 \\
    p19 &    & 17481.6 & 17.08 & 0  & 10 & 10 &    & 50846.5 & 67.12 & 12 & 5  & 17 &    & 83193.2 & 2.15 & 15 & 5  & 20 \\
    p20 &    & 19481.6 & 12.34 & 0  & 10 & 10 &    & 54246.7 & 19.11 & 12 & 5  & 17 &    & 87155.1 & 1.82 & 14 & 5  & 19 \\
    p21 &    & 15199.1 & 77.46 & 0  & 9  & 9  &    & 49458.1 & 66.13 & 15 & 3  & 18 &    & 81630.4 & 1.39 & 17 & 3  & 20 \\
    p22 &    & 14095.4 & 115.28 & 0  & 9  & 9  &    & 46832.6 & 10.31 & 15 & 3  & 18 &    & 78778.4 & 2.04 & 17 & 3  & 20 \\
    p23 &    & 15896 & 57.31 & 0  & 9  & 9  &    & 50311.1 & 16.48 & 14 & 3  & 17 &    & 82778.4 & 1.92 & 17 & 3  & 20 \\
    p24 &    & 17503.6 & 15.59 & 0  & 8  & 8  &    & 53594.6 & 15.77 & 12 & 4  & 16 &    & 86612.4 & 2.48 & 16 & 3  & 19 \\
    p25 &    & 19795.6 & 30.21 & 0  & 9  & 9  &    & 93064.9 & 33.44 & 4  & 8  & 12 &    & 173030 & 12.25 & 5  & 7  & 12 \\
    p26 &    & 18673 & 53.97 & 0  & 10 & 10 &    & 91786.9 & 42.52 & 4  & 8  & 12 &    & 171750 & 15.48 & 5  & 7  & 12 \\
    p27 &    & 20514.4 & 19.14 & 0  & 9  & 9  &    & 94097.3 & 18.68 & 1  & 9  & 10 &    & 174153 & 11.15 & 5  & 7  & 12 \\
    p28 &    & 22314.4 & 9.68 & 0  & 9  & 9  &    & 96100.8 & 6.96 & 2  & 8  & 10 &    & 176551 & 8.31 & 5  & 7  & 12 \\
    p29 &    & 23199.7 & 218.46 & 0  & 14 & 14 &    & 95959 & 68.63 & 1  & 13 & 14 &    & 175722 & 20.70 & 2  & 13 & 15 \\
    p30 &    & 21646.6 & 156.13 & 0  & 14 & 14 &    & 94014.8 & 93.79 & 2  & 13 & 15 &    & 173733 & 31.49 & 2  & 13 & 15 \\
    p31 &    & 24446.6 & 26.78 & 0  & 14 & 14 &    & 96961.9 & 35.57 & 1  & 13 & 14 &    & 176732 & 22.55 & 2  & 13 & 15 \\
    p32 &    & 27077.2 & 15.56 & 0  & 13 & 13 &    & 99726.9 & 12.43 & 0  & 13 & 13 &    & 179716 & 34.24 & 1  & 13 & 14 \\
    p33 &    & 20597.1 & 162.05 & 1  & 10 & 11 &    & 92773.5 & 16.39 & 2  & 10 & 12 &    & 172373 & 10.98 & 2  & 10 & 12 \\
    p34 &    & 19400.5 & 206.43 & 0  & 11 & 11 &    & 91438.2 & 28.34 & 1  & 11 & 12 &    & 171017 & 85.51 & 2  & 10 & 12 \\
    p35 &    & 21582.5 & 199.12 & 0  & 10 & 10 &    & 93838.2 & 18.59 & 1  & 11 & 12 &    & 173416 & 15.10 & 2  & 10 & 12 \\
    p36 &    & 23582.5 & 23.91 & 0  & 10 & 10 &    & 96238.7 & 17.29 & 1  & 11 & 12 &    & 175816 & 7.25 & 2  & 10 & 12 \\
    p37 &    & 18067.1 & 41.02 & 0  & 8  & 8  &    & 91020 & 104.99 & 4  & 6  & 10 &    & 170467 & 47.82 & 6  & 5  & 11 \\
    p38 &    & 17070.1 & 34.33 & 0  & 8  & 8  &    & 89926.4 & 244.87 & 6  & 5  & 11 &    & 169307 & 53.04 & 6  & 5  & 11 \\
    p39 &    & 18670.1 & 45.88 & 0  & 8  & 8  &    & 91807 & 54.75 & 3  & 6  & 9  &    & 171507 & 62.67 & 6  & 5  & 11 \\
    p40 &    & 20252.6 & 41.36 & 0  & 7  & 7  &    & 93607 & 17.05 & 3  & 6  & 9  &    & 173707 & 25.46 & 6  & 5  & 11 \\
    p41 &    & 10641.4 & 2.64 & 0  & 6  & 6  &    & 41371 & 1.24 & 3  & 7  & 10 &    & 73789.5 & 0.63 & 3  & 7  & 10 \\
    p42 &    & 9230.2 & 20.94 & 0  & 5  & 5  &    & 33338.5 & 395.01 & 11 & 4  & 15 &    & 54571.9 & 390.19 & 16 & 3  & 19 \\
    p43 &    & 8058.74 & 33.66 & 0  & 4  & 4  &    & 30081.9 & 10800.00 & 11 & 3  & 14 &    & 48008.9 & 8426.89 & 18 & 3  & 21 \\
    p44 &    & 12065.6 & 0.76 & 0  & 8  & 8  &    & 42859.3 & 0.69 & 3  & 7  & 10 &    & 75596.9 & 0.69 & 3  & 7  & 10 \\
    p45 &    & 10670.3 & 281.57 & 1  & 6  & 7  &    & 33930.2 & 183.68 & 11 & 5  & 16 &    & 53621.2 & 6.59 & 16 & 4  & 20 \\
    p46 &    & 9458.23 & 96.62 & 2  & 4  & 6  &    & 29189 & 10800.00 & 11 & 5  & 16 &    & 44911 & 3925.65 & 20 & 3  & 23 \\
    p47 &    & 11933 & 9.13 & 0  & 9  & 9  &    & 35580.1 & 0.66 & 3  & 7  & 10 &    & 60092.3 & 0.42 & 3  & 7  & 10 \\
    p48 &    & 10584.5 & 1362.60 & 2  & 5  & 7  &    & 25055.4 & 784.10 & 11 & 5  & 16 &    & 37492.1 & 437.16 & 17 & 3  & 20 \\
    p49 &    & 8988.63 & 134.49 & 1  & 5  & 6  &    & 22830.1 & 4708.03 & 10 & 2  & 12 &    & 0.03\% & 10800 & 18 & 2  & 20 \\
    p50 &    & 13831.7 & 1.59 & 0  & 7  & 7  &    & 66864.3 & 0.69 & 2  & 7  & 9  &    & 124561 & 0.57 & 3  & 6  & 9 \\
    p51 &    & 12179.4 & 32.44 & 0  & 7  & 7  &    & 52730.4 & 56.47 & 5  & 7  & 12 &    & 95294.4 & 77.96 & 7  & 6  & 13 \\
    p52 &    & 15983 & 0.61 & 0  & 8  & 8  &    & 74571.7 & 0.26 & 0  & 8  & 8  &    & 138382 & 0.44 & 0  & 8  & 8 \\
    p53 &    & 14715.8 & 7.98 & 0  & 8  & 8  &    & 64795 & 13.68 & 5  & 7  & 12 &    & 117938 & 4.47 & 5  & 7  & 12 \\
    p54 &    & 15502.9 & 0.61 & 0  & 8  & 8  &    & 68880 & 0.49 & 1  & 9  & 10 &    & 126695 & 0.51 & 1  & 9  & 10 \\
    p55 &    & 13571.9 & 6.13 & 0  & 7  & 7  &    & 55608.9 & 14.80 & 1  & 10 & 11 &    & 100424 & 8.90 & 3  & 9  & 12 \\
    \bottomrule
    \end{tabular}%
    }
  \label{tab:impacttransport}%
\end{table}%
\begin{table}[htbp]
  \centering
  \caption{Impact of the production cost for function type 2 on design of the network}
    \resizebox{\textwidth}{!}{
    \begin{tabular}{ccccccccccccccccccc}
    \toprule
       &    & \multicolumn{5}{c}{ $\theta=1$ } &    & \multicolumn{5}{c}{$\theta=10$} &    & \multicolumn{5}{c}{$\theta=20$} \\
\cmidrule{3-7}\cmidrule{9-13}\cmidrule{15-19}    Data  &    & Obj/Gap (\%) & Time & $n_e$ & $n_d$ & $n_T$ &    & Obj/Gap (\%) & Time & $n_e$ & $n_d$ & $n_T$ &    & Obj/Gap (\%) & Time & $n_e$ & $n_d$ & $n_T$ \\
\midrule
   p1 &    & 14468.3 & 0.82 & 0  & 8  & 8  &    & 55789.10 & 4.45 & 0  & 7  & 7  &    & 100531.00 & 3.21 & 0  & 7  & 7 \\
    p2 &    & 13494.2 & 9.48 & 0  & 9  & 9  &    & 54748.60 & 4.54 & 0  & 7  & 7  &    & 99488.60 & 3.76 & 0  & 7  & 7 \\
    p3 &    & 15159.3 & 1.10 & 0  & 8  & 8  &    & 56149.50 & 3.28 & 0  & 7  & 7  &    & 100889.00 & 3.76 & 0  & 7  & 7 \\
    p4 &    & 16599.2 & 0.29 & 0  & 7  & 7  &    & 57547.80 & 2.40 & 0  & 7  & 7  &    & 102288.00 & 2.89 & 0  & 7  & 7 \\
    p5 &    & 15338 & 0.39 & 0  & 10 & 10 &    & 66526.60 & 0.53 & 0  & 10 & 10 &    & 123300.00 & 1.04 & 0  & 10 & 10 \\
    p6 &    & 14024.1 & 0.42 & 0  & 10 & 10 &    & 65212.60 & 0.52 & 0  & 10 & 10 &    & 121984.00 & 0.50 & 0  & 10 & 10 \\
    p7 &    & 16024.1 & 0.43 & 0  & 10 & 10 &    & 67212.60 & 0.49 & 0  & 10 & 10 &    & 123984.00 & 0.47 & 0  & 10 & 10 \\
    p8 &    & 18024 & 0.44 & 0  & 10 & 10 &    & 69212.60 & 0.49 & 0  & 10 & 10 &    & 125985.00 & 0.97 & 0  & 10 & 10 \\
    p9 &    & 13232.2 & 5.04 & 1  & 7  & 8  &    & 53373.60 & 8.88 & 0  & 8  & 8  &    & 97890.20 & 8.31 & 0  & 8  & 8 \\
    p10 &    & 12282 & 14.37 & 1  & 7  & 8  &    & 52423.80 & 9.93 & 0  & 8  & 8  &    & 96940.20 & 8.07 & 0  & 8  & 8 \\
    p11 &    & 13852.1 & 9.98 & 0  & 7  & 7  &    & 54023.60 & 8.03 & 0  & 8  & 8  &    & 98540.20 & 8.16 & 0  & 8  & 8 \\
    p12 &    & 15252.1 & 0.92 & 0  & 7  & 7  &    & 55623.60 & 8.56 & 0  & 8  & 8  &    & 100140.00 & 8.43 & 0  & 8  & 8 \\
    p13 &    & 16876.5 & 72.56 & 0  & 10 & 10 &    & 82191.70 & 147.31 & 0  & 10 & 10 &    & 153509.00 & 136.08 & 0  & 10 & 10 \\
    p14 &    & 15287.7 & 370.08 & 0  & 10 & 10 &    & 79999.80 & 100.07 & 0  & 10 & 10 &    & 151317.00 & 150.09 & 0  & 10 & 10 \\
    p15 &    & 17225.3 & 11.25 & 0  & 9  & 9  &    & 81999.70 & 82.15 & 0  & 10 & 10 &    & 153317.00 & 121.26 & 0  & 10 & 10 \\
    p16 &    & 19025.2 & 2.63 & 0  & 9  & 9  &    & 84000.00 & 83.61 & 0  & 10 & 10 &    & 155317.00 & 208.34 & 0  & 10 & 10 \\
    p17 &    & 16710 & 166.00 & 0  & 10 & 10 &    & 84725.40 & 876.43 & 0  & 11 & 11 &    & 160145.00 & 2286.18 & 0  & 11 & 11 \\
    p18 &    & 15377.4 & 1037.00 & 0  & 11 & 11 &    & 83309.60 & 260.86 & 0  & 11 & 11 &    & 158730.00 & 637.58 & 0  & 11 & 11 \\
    p19 &    & 17481.6 & 17.08 & 0  & 10 & 10 &    & 85509.80 & 343.42 & 0  & 11 & 11 &    & 160929.00 & 1017.26 & 0  & 11 & 11 \\
    p20 &    & 19481.6 & 12.34 & 0  & 10 & 10 &    & 87709.40 & 250.05 & 0  & 11 & 11 &    & 163130.00 & 529.98 & 0  & 11 & 11 \\
    p21 &    & 15199.1 & 77.46 & 0  & 9  & 9  &    & 75938.10 & 343.24 & 0  & 9  & 9  &    & 143337.00 & 673.46 & 0  & 9  & 9 \\
    p22 &    & 14095.4 & 115.28 & 0  & 9  & 9  &    & 74864.90 & 501.23 & 0  & 9  & 9  &    & 142267.00 & 1018.12 & 0  & 9  & 9 \\
    p23 &    & 15896 & 57.31 & 0  & 9  & 9  &    & 76664.90 & 794.89 & 0  & 9  & 9  &    & 144067.00 & 879.28 & 0  & 9  & 9 \\
    p24 &    & 17503.6 & 15.59 & 0  & 8  & 8  &    & 78465.10 & 527.09 & 0  & 9  & 9  &    & 145868.00 & 840.99 & 0  & 9  & 9 \\
    p25 &    & 19795.6 & 30.21 & 0  & 9  & 9  &    & 83001.10 & 136.88 & 0  & 9  & 9  &    & 151183.00 & 214.69 & 0  & 9  & 9 \\
    p26 &    & 18673 & 53.97 & 0  & 10 & 10 &    & 81661.00 & 404.33 & 0  & 9  & 9  &    & 149693.00 & 167.71 & 0  & 9  & 9 \\
    p27 &    & 20514.4 & 19.14 & 0  & 9  & 9  &    & 83461.00 & 604.78 & 0  & 9  & 9  &    & 151493.00 & 234.47 & 0  & 9  & 9 \\
    p28 &    & 22314.4 & 9.68 & 0  & 9  & 9  &    & 85261.20 & 183.07 & 0  & 9  & 9  &    & 153293.00 & 89.07 & 0  & 9  & 9 \\
    p29 &    & 23199.7 & 218.46 & 0  & 14 & 14 &    & 105110.00 & 7153.49 & 0  & 16 & 16 &    & 2.99$\%$ & 10800.00 & 0  & 16 & 16 \\
    p30 &    & 21646.6 & 156.13 & 0  & 14 & 14 &    & 1.34343 per & 10800.00 & 0  & 15 & 15 &    & 2.62$\%$ & 10800.00 & 0  & 16 & 16 \\
    p31 &    & 24446.6 & 26.78 & 0  & 14 & 14 &    & 106460.00 & 9853.85 & 0  & 15 & 15 &    & 2.05$\%$ & 10800.00 & 0  & 15 & 15 \\
    p32 &    & 27077.2 & 15.56 & 0  & 13 & 13 &    & 109460.00 & 6958.09 & 0  & 15 & 15 &    & 2.47$\%$ & 10800.00 & 0  & 15 & 15 \\
    p33 &    & 20597.1 & 162.05 & 1  & 10 & 11 &    & 91558.5 & 10800.00 & 0  & 12 & 12 &    & 2.05$\%$ & 10800.00 & 0  & 12 & 12 \\
    p34 &    & 19400.5 & 206.43 & 0  & 11 & 11 &    & 90129.90 & 5417.24 & 0  & 12 & 12 &    & 2.08$\%$ & 10800.00 & 0  & 12 & 12 \\
    p35 &    & 21582.5 & 199.12 & 0  & 10 & 10 &    & 92529.9 & 4884.20 & 0  & 12 & 12 &    & 3.89$\%$ & 10800.00 & 1  & 11 & 12 \\
    p36 &    & 23582.5 & 23.91 & 0  & 10 & 10 &    & 94930.10 & 3797.02 & 0  & 12 & 12 &    & 2.15$\%$ & 10800.00 & 0  & 12 & 12 \\
    p37 &    & 18067.1 & 41.02 & 0  & 8  & 8  &    & 75738.2 & 4245.00 & 0  & 8  & 8  &    & 2.05$\%$ & 10800.00 & 0  & 8  & 8 \\
    p38 &    & 17070.1 & 34.33 & 0  & 8  & 8  &    & 74741.2 & 5255.92 & 0  & 8  & 8  &    & 2.23$\%$ & 10800.00 & 0  & 8  & 8 \\
    p39 &    & 18670.1 & 45.88 & 0  & 8  & 8  &    & 76341.2 & 2671.45 & 0  & 8  & 8  &    & 2.07$\%$ & 10800.00 & 0  & 8  & 8 \\
    p40 &    & 20252.6 & 41.36 & 0  & 7  & 7  &    & 77941.2 & 4318.77 & 0  & 8  & 8  &    & 2.04$\%$ & 10800.00 & 0  & 8  & 8 \\
    p41 &    & 10641.4 & 2.64 & 0  & 6  & 6  &    & 41918.2 & 6.41 & 0  & 5  & 5  &    & 76583.00 & 4.21 & 0  & 5  & 5 \\
    p42 &    & 9230.2 & 20.94 & 0  & 5  & 5  &    & 36169.2 & 38.44 & 0  & 4  & 4  &    & 65055.30 & 22.13 & 0  & 4  & 4 \\
    p43 &    & 8058.74 & 33.66 & 0  & 4  & 4  &    & 32355.1 & 121.57 & 0  & 4  & 4  &    & 59250.50 & 557.56 & 0  & 4  & 4 \\
    p44 &    & 12065.6 & 0.76 & 0  & 8  & 8  &    & 52359.8 & 4.08 & 0  & 7  & 7  &    & 96405.70 & 2.67 & 0  & 6  & 6 \\
    p45 &    & 10670.3 & 281.57 & 1  & 6  & 7  &    & 43645.4 & 46.06 & 0  & 5  & 5  &    & 78962.20 & 27.62 & 0  & 5  & 5 \\
    p46 &    & 9458.23 & 96.62 & 2  & 4  & 6  &    & 39617.50 & 1268.24 & 1  & 4  & 5  &    & 3.12$\%$ & 10800.70 & 1  & 4  & 5 \\
    p47 &    & 11933 & 9.13 & 0  & 9  & 9  &    & 52069.7 & 3.57 & 0  & 7  & 7  &    & 96031.10 & 6.08 & 0  & 6  & 6 \\
    p48 &    & 10584.5 & 1362.60 & 2  & 5  & 7  &    & 44592.4 & 44.26 & 0  & 5  & 5  &    & 81142.50 & 28.66 & 0  & 5  & 5 \\
    p49 &    & 8988.63 & 134.49 & 1  & 5  & 6  &    & 37674.9 & 476.34 & 0  & 4  & 4  &    & 68696.50 & 435.48 & 0  & 4  & 4 \\
    p50 &    & 13831.7 & 1.59 & 0  & 7  & 7  &    & 50842.1 & 4.02 & 0  & 6  & 6  &    & 90784.10 & 4.31 & 0  & 6  & 6 \\
    p51 &    & 12179.4 & 32.44 & 0  & 7  & 7  &    & 47847 & 329.77 & 0  & 6  & 6  &    & 84839.40 & 27.12 & 0  & 5  & 5 \\
    p52 &    & 15983 & 0.61 & 0  & 8  & 8  &    & 62944.2 & 6.65 & 0  & 8  & 8  &    & 113596.00 & 2.69 & 0  & 8  & 8 \\
    p53 &    & 14715.8 & 7.98 & 0  & 8  & 8  &    & 57434.9 & 38.12 & 0  & 6  & 6  &    & 101736.00 & 27.25 & 0  & 6  & 6 \\
    p54 &    & 15502.9 & 0.61 & 0  & 8  & 8  &    & 64603.6 & 1.81 & 0  & 8  & 8  &    & 118979.00 & 1.05 & 0  & 8  & 8 \\
    p55 &    & 13571.9 & 6.13 & 0  & 7  & 7  &    & 57276.1 & 21.16 & 0  & 7  & 7  &    & 105286.00 & 26.72 & 0  & 7  & 7 \\
\bottomrule
    \end{tabular}%
    }
  \label{tab:impactprod}%
\end{table}%
\begin{table}[htbp]
  \centering
  \caption{Impact of the economic point location for function type 2 on the design of the network}
     \resizebox{\textwidth}{!}{
    \begin{tabular}{ccccccccccccccccccc}
    \toprule
       &    & \multicolumn{5}{c}{ $\beta=0.5$ } &    & \multicolumn{5}{c}{$\beta=0.6$} &    & \multicolumn{5}{c}{$\beta =0.75$} \\
\cmidrule{3-7}\cmidrule{9-13}\cmidrule{15-19}    Data  &    & Obj/Gap (\%) & Time & $n_e$ & $n_d$ & $n_T$ &    & \multicolumn{1}{c}{Obj/Gap (\%)} & \multicolumn{1}{c}{Time} & \multicolumn{1}{c}{$n_e$} & \multicolumn{1}{c}{$n_d$} & \multicolumn{1}{c}{$n_T$} &    & \multicolumn{1}{c}{Obj/Gap (\%)} & \multicolumn{1}{c}{Time} & \multicolumn{1}{c}{$n_e$} & \multicolumn{1}{c}{$n_d$} & \multicolumn{1}{c}{$n_T$} \\
\cmidrule{1-1}\cmidrule{3-7}\cmidrule{9-13}\cmidrule{15-19}    p1 &    & 14468.3 & 0.82 & 0  & 8  & 8  &    & 14232.7 & 1.43308 & 0  & 8  & 8  &    & 13680.1 & 30.2327 & 4  & 4  & 8 \\
    p2 &    & 13494.2 & 9.48 & 0  & 9  & 9  &    & 13311.1 & 16.4279 & 3  & 6  & 9  &    & 12730.6 & 57.3143 & 0  & 8  & 8 \\
    p3 &    & 15159.3 & 1.10 & 0  & 8  & 8  &    & 14915.9 & 2.01412 & 0  & 7  & 7  &    & 14207.6 & 14.0928 & 2  & 5  & 7 \\
    p4 &    & 16599.2 & 0.29 & 0  & 7  & 7  &    & 16315.9 & 0.554033 & 0  & 7  & 7  &    & 15607.6 & 3.4812 & 0  & 7  & 7 \\
    p5 &    & 15338 & 0.39 & 0  & 10 & 10 &    & 15061.4 & 0.667041 & 0  & 10 & 10 &    & 14435.6 & 84.6653 & 6  & 4  & 10 \\
    p6 &    & 14024.1 & 0.42 & 0  & 10 & 10 &    & 13747.4 & 0.640039 & 0  & 10 & 10 &    & 13121.6 & 107.495 & 5  & 5  & 10 \\
    p7 &    & 16024.1 & 0.43 & 0  & 10 & 10 &    & 15747.4 & 0.847054 & 0  & 10 & 10 &    & 15050.7 & 16.565 & 0  & 9  & 9 \\
    p8 &    & 18024 & 0.44 & 0  & 10 & 10 &    & 17747.4 & 0.814053 & 0  & 10 & 10 &    & 16850.7 & 4.94229 & 0  & 9  & 9 \\
    p9 &    & 13232.2 & 5.04 & 1  & 7  & 8  &    & 13112.6 & 46.8458 & 2  & 6  & 8  &    & 12678 & 25.0444 & 5  & 2  & 7 \\
    p10 &    & 12282 & 14.37 & 1  & 7  & 8  &    & 12162.7 & 38.6786 & 4  & 4  & 8  &    & 11865.1 & 35.176 & 6  & 1  & 7 \\
    p11 &    & 13852.1 & 9.98 & 0  & 7  & 7  &    & 13654.5 & 9.71325 & 0  & 7  & 7  &    & 13234.7 & 63.0236 & 0  & 6  & 6 \\
    p12 &    & 15252.1 & 0.92 & 0  & 7  & 7  &    & 15054.5 & 4.34199 & 1  & 6  & 7  &    & 14435.4 & 3.83422 & 0  & 6  & 6 \\
    p13 &    & 16876.5 & 72.56 & 0  & 10 & 10 &    & 16512.8 & 94.1009 & 0  & 10 & 10 &    & 15515.1 & 79.9027 & 2  & 7  & 9 \\
    p14 &    & 15287.7 & 370.08 & 0  & 10 & 10 &    & 14998.1 & 1113.59 & 0  & 10 & 10 &    & 14098.7 & 2727.35 & 3  & 6  & 9 \\
    p15 &    & 17225.3 & 11.25 & 0  & 9  & 9  &    & 16862.8 & 24.7856 & 0  & 9  & 9  &    & 15898.7 & 448.783 & 1  & 8  & 9 \\
    p16 &    & 19025.2 & 2.63 & 0  & 9  & 9  &    & 18662.8 & 7.48348 & 0  & 9  & 9  &    & 17698.7 & 402.31 & 2  & 7  & 9 \\
    p17 &    & 16710 & 166.00 & 0  & 10 & 10 &    & 16354.2 & 314.149 & 0  & 10 & 10 &    & 15414.4 & 200.135 & 0  & 9  & 9 \\
    p18 &    & 15377.4 & 1037.00 & 0  & 11 & 11 &    & 0.36$\%$ & 10800 & 2  & 9  & 11 &    & 1.09$\%$ & 10800.3 & 5  & 5  & 10 \\
    p19 &    & 17481.6 & 17.08 & 0  & 10 & 10 &    & 17132.4 & 26.5009 & 0  & 10 & 10 &    & 16131.4 & 45.0325 & 0  & 9  & 9 \\
    p20 &    & 19481.6 & 12.34 & 0  & 10 & 10 &    & 19132.5 & 22.3376 & 0  & 10 & 10 &    & 17931.7 & 37.2431 & 0  & 9  & 9 \\
    p21 &    & 15199.1 & 77.46 & 0  & 9  & 9  &    & 14946.1 & 362.867 & 0  & 8  & 8  &    & 14257.4 & 717.748 & 5  & 3  & 8 \\
    p22 &    & 14095.4 & 115.28 & 0  & 9  & 9  &    & 13887.1 & 3107.53 & 1  & 8  & 9  &    & 13305.7 & 9143.93 & 4  & 4  & 8 \\
    p23 &    & 15896 & 57.31 & 0  & 9  & 9  &    & 15592.4 & 22.1514 & 0  & 8  & 8  &    & 14905.8 & 1952.36 & 6  & 2  & 8 \\
    p24 &    & 17503.6 & 15.59 & 0  & 8  & 8  &    & 17192.4 & 14.77 & 0  & 8  & 8  &    & 16505.7 & 1991.06 & 5  & 3  & 8 \\
    p25 &    & 19795.6 & 30.21 & 0  & 9  & 9  &    & 19478.6 & 37.0744 & 0  & 9  & 9  &    & 18762.7 & 5324.71 & 3  & 6  & 9 \\
    p26 &    & 18673 & 53.97 & 0  & 10 & 10 &    & 18397.8 & 228.605 & 0  & 9  & 9  &    & 17681.8 & 2065.93 & 4  & 5  & 9 \\
    p27 &    & 20514.4 & 19.14 & 0  & 9  & 9  &    & 20197.8 & 33.9302 & 0  & 9  & 9  &    & 19460 & 1493.36 & 6  & 2  & 8 \\
    p28 &    & 22314.4 & 9.68 & 0  & 9  & 9  &    & 21960.8 & 41.9831 & 0  & 8  & 8  &    & 21060 & 1043.66 & 3  & 5  & 8 \\
    p29 &    & 23199.7 & 218.46 & 0  & 14 & 14 &    & 22782 & 2077.35 & 0  & 13 & 13 &    & 21434.8 & 10801.8 & 5  & 8  & 13 \\
    p30 &    & 21646.6 & 156.13 & 0  & 14 & 14 &    & 21245.5 & 352.243 & 0  & 14 & 14 &    & 20002 & 342.178 & 2  & 11 & 13 \\
    p31 &    & 24446.6 & 26.78 & 0  & 14 & 14 &    & 23954.7 & 90.4671 & 0  & 13 & 13 &    & 22602 & 472.166 & 0  & 13 & 13 \\
    p32 &    & 27077.2 & 15.56 & 0  & 13 & 13 &    & 26554.7 & 56.6318 & 0  & 13 & 13 &    & 25202 & 1472.61 & 2  & 11 & 13 \\
    p33 &    & 20597.1 & 162.05 & 1  & 10 & 11 &    & 20281 & 1489.02 & 1  & 10 & 11 &    & 0.01$\%$ & 10804.9 & 4  & 6  & 10 \\
    p34 &    & 19400.5 & 206.43 & 0  & 11 & 11 &    & 19099.3 & 2740.54 & 0  & 11 & 11 &    & 0.18$\%$ & 10801.6 & 1  & 9  & 10 \\
    p35 &    & 21582.5 & 199.12 & 0  & 10 & 10 &    & 21168.5 & 129.446 & 0  & 10 & 10 &    & 0.15$\%$ & 10800.4 & 1  & 9  & 10 \\
    p36 &    & 23582.5 & 23.91 & 0  & 10 & 10 &    & 23168.5 & 152.496 & 0  & 10 & 10 &    & 0.09$\%$ & 10828.7 & 0  & 9  & 9 \\
    p37 &    & 18067.1 & 41.02 & 0  & 8  & 8  &    & 17887.1 & 435.53 & 1  & 7  & 8  &    & 17340.6 & 3744.78 & 3  & 4  & 7 \\
    p38 &    & 17070.1 & 34.33 & 0  & 8  & 8  &    & 16890.5 & 329.547 & 2  & 6  & 8  &    & 16470.6 & 3949.67 & 4  & 3  & 7 \\
    p39 &    & 18670.1 & 45.88 & 0  & 8  & 8  &    & 18490.5 & 511.143 & 1  & 7  & 8  &    & 17870.6 & 2794.19 & 3  & 4  & 7 \\
    p40 &    & 20252.6 & 41.36 & 0  & 7  & 7  &    & 19961 & 137.124 & 0  & 7  & 7  &    & 19270.8 & 825.882 & 5  & 2  & 7 \\
    p41 &    & 10641.4 & 2.64 & 0  & 6  & 6  &    & 10549.9 & 62.3087 & 2  & 4  & 6  &    & 10170.8 & 12.5898 & 0  & 5  & 5 \\
    p42 &    & 9230.2 & 20.94 & 0  & 5  & 5  &    & 9147.41 & 264.25 & 3  & 2  & 5  &    & 8851.97 & 202.35 & 0  & 4  & 4 \\
    p43 &    & 8058.74 & 33.66 & 0  & 4  & 4  &    & 7973.14 & 78.0705 & 1  & 3  & 4  &    & 7861.35 & 138.409 & 3  & 1  & 4 \\
    p44 &    & 12065.6 & 0.76 & 0  & 8  & 8  &    & 11895.3 & 2.39914 & 0  & 8  & 8  &    & 11532.7 & 147.139 & 6  & 1  & 7 \\
    p45 &    & 10670.3 & 281.57 & 1  & 6  & 7  &    & 10565.9 & 474.507 & 1  & 5  & 6  &    & 10252.2 & 1693.56 & 1  & 5  & 6 \\
    p46 &    & 9458.23 & 96.62 & 2  & 4  & 6  &    & 9367.1 & 187.329 & 2  & 4  & 6  &    & 9225.52 & 869.016 & 6  & 0  & 6 \\
    p47 &    & 11933 & 9.13 & 0  & 9  & 9  &    & 11781.6 & 64.1529 & 1  & 8  & 9  &    & 11372.6 & 55.8354 & 3  & 5  & 8 \\
    p48 &    & 10584.5 & 1362.60 & 2  & 5  & 7  &    & 0.024$\%$ & 10800.7 & 3  & 5  & 8  &    & 0.012$\%$ & 10800.5 & 4  & 2  & 6 \\
    p49 &    & 8988.63 & 134.49 & 1  & 5  & 6  &    & 8915.87 & 807.551 & 2  & 4  & 6  &    & 8723.42 & 1083.92 & 4  & 1  & 5 \\
    p50 &    & 13831.7 & 1.59 & 0  & 7  & 7  &    & 13624.8 & 2.06315 & 0  & 7  & 7  &    & 13128.8 & 11.5238 & 4  & 3  & 7 \\
    p51 &    & 12179.4 & 32.44 & 0  & 7  & 7  &    & 12002.4 & 52.1286 & 0  & 6  & 6  &    & 11503.7 & 360.451 & 5  & 1  & 6 \\
    p52 &    & 15983 & 0.61 & 0  & 8  & 8  &    & 15626.9 & 0.549042 & 0  & 8  & 8  &    & 14717.2 & 0.866063 & 0  & 8  & 8 \\
    p53 &    & 14715.8 & 7.98 & 0  & 8  & 8  &    & 14490.3 & 37.5556 & 0  & 8  & 8  &    & 13723.8 & 17.9553 & 0  & 7  & 7 \\
    p54 &    & 15502.9 & 0.61 & 0  & 8  & 8  &    & 15240.9 & 1.49411 & 0  & 8  & 8  &    & 14556.5 & 17.8293 & 3  & 5  & 8 \\
    p55 &    & 13571.9 & 6.13 & 0  & 7  & 7  &    & 13344.4 & 11.7128 & 0  & 7  & 7  &    & 12716.7 & 50.3542 & 0  & 7  & 7 \\
    \bottomrule
    \end{tabular}%
    }
  \label{tab:impact_eco_p}%
\end{table}%

\end{document}